\documentclass[ twoside, 11pt]{article}
\usepackage[margin=1in]{geometry}

\usepackage{setspace}
\usepackage{fancyhdr}

\newcommand\shorttitle{Sums of m-dependent
random variables}
\newcommand\authors{
V.~\v{C}ekanavi\v{c}ius, P. Vellaisamy,}

\usepackage{amsmath}
\usepackage{amssymb}
\usepackage[square]{natbib}
\setcitestyle{numbers}

\newtheorem{theorem}{Theorem}[section]
\newtheorem{corollary}{Corollary}[section]
\newtheorem{lemma}{Lemma}[section]

\newcommand{\Proof}{\textbf{Proof. }}            
\newcommand{\qed}{$\square$}                     

\newcommand{\RR}{\mathbb{R}}

\newcommand{\ZZ}{\mathbb{Z}}
\newcommand{\Zpl}{\ZZ_+}

\newcommand{\ii}{{\mathrm i}}
\newcommand{\ee}{{\mathrm e}}

\newcommand{\dd}{{\mathrm{d}}}

\newcommand{\dirac}{I}               

\newcommand{\BI}{\mathrm{Bi}}    
\newcommand{\NB}{\mathrm{NB}}   
\newcommand{\Pois}{\mathrm{Pois}}    
\newcommand{\G}{\mathrm{G}}    

\newcommand{\Ha}{\mathrm{H}}    

\newcommand{\norm}[1]{\|#1\|}                    
\newcommand{\Norm}[1]{\Big\|#1\Big\|}            
\newcommand{\ab}[1]{\vert#1\vert}                
\newcommand{\Ab}[1]{\Big\vert#1\Big\vert}        

\newcommand{\exponent}[1]{\exp\{#1\}}            
\newcommand{\Exponent}[1]{\exp\Bigl\{#1\Bigr\}}  

\newcommand{\eit}{\ee^{\ii t}}                   

\newcommand{\Expect}{\mathrm{E}}                 
\newcommand{\Var}{\mathrm{Var}}                  

\newcommand{\qubar}{\overline{q}}
\newcommand{\pbar}{\overline{p}}

\newcommand{\vfi}{\varphi}
\newcommand{\wE}{\w\Expect}
\newcommand{\wexp}{\w\Expect^{+}}
\newcommand{\eL}{{\cal L}}

\newcommand{\w}{\widehat}

\newcommand{\floor}[1]{{\lfloor #1\rfloor}}

\newcommand{\binomial}[2]{\genfrac{(}{)}{0pt}{}{#1}{#2}}

\begin{document}
\title{Discrete approximations for sums of m-dependent
random variables}

\author{V. \v Cekanavi\v cius  and P. Vellaisamy      \\
{\small
Department of Mathematics and Informatics, Vilnius University,}\\
{\small Naugarduko 24, Vilnius 03225, Lithuania.}\\{\small E-mail:
vydas.cekanavicius@mif.vu.lt } \\{\small and}
\\{\small
 Department of Mathematics, Indian Institute of Technology Bombay,} \\
 {\small Powai, Mumbai-
400076, India.}\\{\small  E-mail: pv@math.iitb.ac.in} }
\date{}

\maketitle

\begin{abstract}

Sums of  $m$-dependent integer-valued random variables are
approximated by compound Poisson, negative binomial and binomial
distributions and signed compound Poisson measures. Estimates are
obtained for the total variation  metric. The results are
then applied to statistics of $m$-dependent $(k_1,k_2)$ events and
2-runs. Heinrich's method and smoothing properties of convolutions
are used for the proofs.

\vspace*{.5cm} \noindent {\emph{Key words:} \small  Compound
Poisson distribution, signed compound Poisson measure, negative
binomial, binomial, m-dependent variables, total variation norm.}

\vspace*{.5cm} \noindent {\small {\it MSC 2000 Subject
Classification}:
Primary 60F05.   
Secondary 60G50;     
}
\end{abstract}

\newpage

\section{The setup} 

In this paper, we consider  sums $S_n= X_1+X_2+\cdots+X_n$ of
non-identically distributed 1-dependent random variables
concentrated on nonnegative integers. Our aim is to estimate the
closeness of $S_n$ to  compound Poisson, negative binomial and
binomial distributions,  under some conditions
for factorial moments. For the proof of the main results, we use
Heinrich's \citep{H82},\citep{H87} version of the characteristic
function method. Though this method does not allow to obtain small
absolute constants, it is flexible enough for obtaining
 asymptotically sharp constants, as demonstrated for 2-runs statistic. Moreover, our approach allows
for construction of asymptotic expansions.

 We recall that the sequence of random variables
$\{X_k \}_{k \geq 1}$
  is called $m$-dependent if, for $1 < s < t < \infty$, $t- s > m$, the sigma algebras
generated by $X_1,\dots,X_s$ and $X_t,X_{t+1}\dots$ are
independent. It is clear that,  by grouping consecutive summands, we
can reduce the sum of $m$-dependent variables to the sum of
1-dependent ones. Therefore,  the results of this paper can  be
applied  for some cases of $m$-dependent variables, as exemplified
by binomial approximation to $(k_1,k_2)$ events.

Let us introduce some necessary notations. Let $\{Y_{k}\}_{k \geq 1}$
be a sequence of arbitrary real or complex-valued random
variables. We assume that $\w\Expect (Y_1)=\Expect Y_1$ and, for
$k\geqslant 2$, define $\w\Expect (Y_1,Y_2,\cdots Y_k)$  by
\begin{equation*}
 \w\Expect (Y_1,Y_2,\cdots, Y_k)=\Expect Y_1Y_2\cdots
Y_k-\sum_{j=1}^{k-1}\w\Expect (Y_1,\cdots ,Y_j)\Expect
Y_{j+1}\cdots Y_{k}. \label{capY}
\end{equation*}

We define $j$-th factorial moment of $X_k$ by $\nu_j(k)=\Expect
X_k(X_k-1)\cdots(X_k-j+1)$, ($k=1,2,\dots,n$, $j=1,2,\dots$).
Let
\begin{eqnarray*} \Gamma_1&=&\Expect
S_n=\sum_{k=1}^n\nu_1(k),\quad
\Gamma_2=\frac{1}{2}(\Var S_n-\Expect S_n)=\frac{1}{2}\sum_{k=1}^n\big(\nu_2(k)-\nu_1^2(k)\big)+\sum_{k=2}^n\wE(X_{k-1},X_k),\\
\Gamma_3&=&\frac{1}{6}\sum_{k=1}^n\big(\nu_3(k)-3\nu_1(k)\nu_2(k)+2\nu_1^3(k)\big)
-\sum_{k=2}^n\big(\nu_1(k-1)+\nu_1(k)\big)\wE(X_{k-1},X_k)\\
&& +\frac{1}{2}\sum_{k=2}^n\big(\wE(X_{k-1}(X_{k-1}-1),X_k)
+\wE(X_{k-1},X_k(X_k-1))\big)
+\sum_{k=3}^n\wE(X_{k-2},X_{k-1},X_k).
\end{eqnarray*}

For the sake of convenience, we assume that $X_k\equiv 0$ and
$\nu_j(k)=0$ if $k\leqslant 0$ and $\sum_k^n=0$ if $k>n$. We
denote the distribution and characteristic function of $S_n$ by
$F_n$ and $\w F_{n}(t)$, respectively. Below we show that $\Gamma_1$,
$2\Gamma_2$ and $6\Gamma_3$ are factorial cumulants of $F_n$, that is,
\begin{equation*}\label{gamma}
 \w F_n(t)=\exponent{\Gamma_1(\eit-1)+\Gamma_2(\eit-1)^2+\Gamma_3(\eit-1)^3+\dots}.
\end{equation*}
For approximation of $F_n$, it is natural to use measures or
distributions which allow similar expressions.

 Let $\dirac_a$ denote the distribution concentrated at real $a$
and set $\dirac=\dirac_0$. Henceforth, the products and powers of measures are
understood in the convolution sense. Further, for a measure $M$, we
set $M^0=\dirac$ and
\[\ee^M:=\exponent{M}=\sum_{k=0}^\infty\frac{1}{k!}\,M^k.\]
 The total variation norm  of measure $M$ is denoted by
\[\norm{M}=\sum_{k=-\infty}^\infty\ab{M\{k\}}.\]
 We use symbol $C$ to
denote all (in general, different) positive absolute constants.
We use symbols $\theta$ and
$\Theta$ to denote all real or complex quantities satisfying
$\ab{\theta}\leqslant 1$ and all measures of finite variation
satisfying $\norm{\Theta}=1$, respectively.

Next we define  approximations of this
paper. Let
\[
\Pois(\Gamma_1)=\exponent{\Gamma_1(\dirac_1-\dirac)},\quad
 \G=\exponent{\Gamma_1(\dirac_1-\dirac)+\Gamma_2(\dirac_1-\dirac)^2}.
 \]
 It is easy to see that $\Pois(\Gamma_1)$ is Poisson
distribution with  parameter $\Gamma_1$. In general, $\G$ is a
signed measure, since $\Gamma_2$ can be negative. Signed compound
Poisson measures similar to $\G$ are used in numerous papers, see
\citep{BC02}, \citep{BaXi99}, \citep{CV10},
\citep{Roo03},  and the references therein. In comparison to the
Poisson distribution, the main benefit of $\G$   is
matching of two moments, which then allows for the accuracy
comparable to the one achieved by the normal approximation.
 This fact is illustrated in
the next two sections. From a practical point of view, signed
measures are not always convenient to use, since for calculation
of their 'probabilities' one needs  inverse Fourier transform or
recursive algorithms. Therefore, we also prove estimates for such
widely used distributions as binomial and negative binomial. We define
the binomial distribution of this paper as
\begin{equation*}\label{binom}
\BI(N,\bar{p})=(\dirac+\bar{p}(\dirac_1-\dirac))^N,\quad
N=\floor{\tilde N}, \quad \tilde N
=\frac{\Gamma_1^2}{2\ab{\Gamma_2}},
\quad\bar{p}=\frac{\Gamma_1}{N}.
\end{equation*}
 Here, we use
$\floor{\tilde N}$ to denote the integer part of $\tilde N$, that
is, $\tilde N = N+\epsilon$, for some $ 0\leq \epsilon <1$. Also, we define
negative binomial distribution and choose its parameters in the
following way:
\begin{equation}\label{nbinom}
\NB(r,\bar{q})\{j\}=\frac{\Gamma(r+j)}{j!\Gamma(r)}\,\bar{q}^r(1-\bar{q})^j,\quad
(j\in\Zpl),\qquad \frac{r(1-\bar{q})}{\bar{q}}=\Gamma_1,\quad
r\bigg(\frac{1-\bar{q}}{\bar{q}}\bigg)^2=2\Gamma_2.
\end{equation}
Note that  symbols $\bar{q}$ and $\bar{p}$ are not related and, in
general, $\bar{q}+\bar{p}\neq 1$.


\section{Known results}

There are many results dealing with approximations to the sum of
dependent integer-valued random variables. Note, however, that
with very few exceptions: a) all papers are devoted to the sums of
indicator variables only; b) results are not related to
$k$-dependent variables. For example, indicators connected in a
Markov chain are investigated in \citep{CV10}, \citep{XZ09}.  The
most general results, containing $k$-dependent variables as
partial cases, are obtained for birth-death processes with some
stochastic ordering, see \citep{BrX01}, \citep{Daly12},
\citep{Eich99} and the references therein.

 Arguably the best explored case of sums of
$k$-dependent integer-valued random variables is $k$-runs.
Let $\eta_i\sim Be(p_i)$ (i=1,2,\dots) be
independent Bernoulli variables. Let us define
$\xi_i=\prod_{s=i}^{i+k-1}\eta_i$, $S^{*}=\sum_{i=1}^n\xi_i$.
 The sum $S^{*}$ is called $k$-runs statistic. Note that frequently $\eta_{i+nm}$ is treated as $\eta_i$ for $1\leqslant
i\leqslant n$ and $m=\pm 1, \pm 2,\dots$.
Approximations of 2 or  $k$-runs statistic by Poisson,  negative binomial distribution or signed compound Poisson
measure are considered in \citep{BaXi99},
\citep{BrX01},\citep{Daly12}, \citep{Rol05} and \citep{WXia08}.
Particularly in \citep{BrX01} it was proved that,
if $k=2$ and $p_i\equiv p$,  $n\geqslant 2$ and $p<2/3$,
then
\begin{equation}
\norm{\eL(S^{*})-\NB(\tilde r,\tilde q)}\leqslant
\frac{64.4p}{\sqrt{(n-1)(1-p)^3}}. \label{BXNB}
\end{equation}
Here  $\tilde q=(2p-3p^2)/(1+2p-3p^2)$ and
$(1-\tilde{q})/\tilde{q}=np^2$.

The $k$-runs statistic has very explicit dependency of summands.
Meanwhile, our aim is  to obtain a general result which includes sums of
independent random variables  as a particular case. Except for
examples, no specific assumptions about the structure of summands
are made. For bounded and identically distributed random variables
a similar approach is taken in \citep{PeCe2}. We
 give one example from \citep{PeCe2} in the notation of the  previous Section.
 Let the $X_i$ be identically distributed,
$\ab{X_1}\leqslant C$, and, for $n\to\infty$,
\vspace*{-0.2cm}
 \begin{equation}
\nu_1(1)=o(1),\quad\nu_2(1)=o(\nu_1(1)),\quad \Expect
X_1X_2=o(\nu_1(1)), \quad n\nu_1(1)\to\infty. \label{cond}
 \end{equation}
Then
\begin{equation*}
\norm{F_n-G}=O\bigg(\frac{\tilde
R}{\nu_1(1)\sqrt{n\nu_1(1)}}\bigg),\label{PetrC}
\end{equation*}
where
\begin{equation*}
\tilde R =\nu_3(1)+\nu_1(1)\nu_2(1)+\nu_1^3(1)+\Expect
(X_1(X_1-1)X_2+X_1X_2(X_2-1))+\nu_1(1)\Expect X_1X_2+\Expect
X_1X_2X_3.
\end{equation*}

\vspace*{-0.2cm}
 Condition (\ref{cond}) implies that $X_i$ form a triangular array and
  $P(X_i=k)=o(1)$, $k>1$.
  Thus, the classical case of a sequence of random variables, so
  typical for CLT,
 is completely excluded.
 Moreover, assumption $\ab{X_1}\leqslant C$ seems rather
strong. For example, then one can not consider  Poisson or
geometric random variables as possible summands.


\section{Results}

 All results are obtained under the following conditions:
\begin{eqnarray}
\nu_1(k)&\leqslant& 1/100, \quad
\nu_2(k)\leqslant \nu_1(k),\quad \nu_4(k)<\infty,\quad (k=1,2,\dots,n),\label{nu12}\\
\lambda &:=&
\sum_{k=1}^n\nu_1(k)-1.52\sum_{k=1}^n\nu_2(k)-12\sum_{k=2}^n\Expect
X_{k-1}X_{k}>0.\label{lambda}
\end{eqnarray}
The last condition is satisfied, if the following two assumptions
hold
 \begin{equation}
\sum_{k=1}^n\nu_2(k)\leqslant\frac{\Gamma_1}{20},\qquad
\sum_{k=2}^n\ab{Cov(X_{k-1},X_k)}\leqslant\frac{\Gamma_1}{20}.\label{3ab}
 \end{equation}
Moreover, if (\ref{nu12}) and (\ref{3ab}) hold, then $\lambda
> 0.2\Gamma_1$. Indeed, then
\[ \Expect X_{k-1}X_{k}\leqslant \ab{Cov(X_{k-1},X_k)+\nu_1(k-1)\nu_1(k)}\leqslant \ab{Cov(X_{k-1},X_k)}+0.01\nu_1(k). \]

Conditions  above are weaker than
(\ref{cond}). For example, $X_j$ are not necessarily bounded by
some absolute constant.

Next we define  remainder terms. Let
\begin{eqnarray*}
R_0&=&\sum_{k=1}^n\Big\{\nu_2(k)+\nu_1^2(k)+\Expect X_{k-1}X_k\Big\},\\
R_1&=&\sum_{k=1}^n\Big\{\nu_1^3(k)+\nu_1(k)\nu_2(k)+\nu_3(k)+[\nu_1(k-2)+\nu_1(k-1)+\nu_1(k)]\Expect
X_{k-1}X_k\nonumber\\
&&+\wE^{+}_2(X_{k-1},X_k)+\wE^{+}(X_{k-2},X_{k-1},X_k)\Big\},
\\
R_2&=&\sum_{k=1}^n\Big\{\nu_1^4(k)+\nu_2^2(k)+\nu_4(k)+[\nu_1(k-2)+\nu_1(k-1)+\nu_1(k)][\nu_3(k)+\wE^{+}_2(X_{k-1},X_k)]\nonumber\\
&&+\Big(\Expect
X_{k-1}X_k\Big)^2+\sum_{l=0}^3\nu_1(k-l)\wE^{+}(X_{k-2},X_{k-1},X_{k})+\wE^{+}_2(X_{k-2},X_{k-1},X_k)\nonumber\\
&&+\wE^{+}_3(X_{k-1},X_{k})+\wE^{+}(X_{k-3},X_{k-2},X_{k-1},X_k)\Big\}.
\end{eqnarray*}
Here
\begin{eqnarray*}
\wE^{+}(X_1)&=&\Expect X_1,\qquad\wE^{+}(X_1,X_2)=\Expect
X_1X_2+\Expect X_1\Expect X_2,\nonumber\\
 \wE^{+}(X_1,\dots,X_k)&=&\Expect X_1\dots
X_k+\sum_{j=1}^{k-1}\wE^{+}(X_1,X_2,\dots,X_j)\Expect
X_{j+1}X_{j+2}\cdots X_k,\\
\wE^{+}_2(X_{k-1},X_k)&=&\wexp(X_{k-1}(X_{k-1}-1),X_k)+\wexp(X_{k-1},X_k(X_k-1)),\\
\wexp_2(X_{k-2},X_{k-1},X_k)&=&\wexp(X_{k-2}(X_{k-2}-1),X_{k-1},X_k)+\wexp(X_{k-2},X_{k-1}(X_{k-1}-1),X_k)\nonumber\\&&+
\wexp(X_{k-2},X_{k-1},X_k(X_k-1)),\\
\wexp_3(X_{k-1},X_k)&=&\wexp(X_{k-1}(X_{k-1}-1)(X_{k-1}-2),X_k)+
\wexp(X_{k-1}(X_{k-1}-1),X_k(X_{k}-1))\nonumber\\&&+\wexp(X_{k-1},X_k(X_k-1)(X_k-2)).
\end{eqnarray*}

For better understanding of the order of remainder terms, let us
consider the case of Bernoulli variables
$P(X_i=1)=1-P(X_i=0)=p_i$. If all $X_i$ are independent, then
$R_0=C\sum_1^n p_i^2$ and $R_1=C\sum_1^n p_i^{3}$. If $X_i$ are
1-dependent, then at least $R_0\leqslant C\sum_1^n p_i$ and
$R_1\leqslant C\sum_{1}^n p_k^{3/2}$. If some additional
information about $X_i$ is available (for example, that they form
2-runs), then the estimates are somewhat in between.

Our aim is investigation of approximations with at least two
parameters. However, for the sake of completeness, we start with the
Poisson approximation. Note that Poisson approximation (for
indicator variables) is considered in \citep{AGG90}, \citep{BHJ92}
under much more general conditions than  assumed in this paper.

\begin{theorem} \label{PoissonT} Let conditions (\ref{nu12}) and (\ref{lambda}) be
satisfied. Then, for all $n$,
\begin{eqnarray}
\norm{F_n-\Pois(\Gamma_1)}&\leqslant&
C R_0\big\{1+\Gamma_1\min(1,\lambda^{-1})\big\}\min(1,\lambda^{-1}),\label{PV}
\\
\norm{F_n-\Pois(\Gamma_1)(\dirac+\Gamma_2(\dirac_1-\dirac)^2)}&\leqslant&
C \big\{1+\Gamma_1\min(1,\lambda^{-1})\big\}(R_0^2\min(1,\lambda^{-2})\nonumber\\&&+R_1\min(1,\lambda^{-3/2})).\label{PVA}
\end{eqnarray}
\end{theorem}

If all $X_i\sim Be(1,p_i)$  are independent, then the order of
accuracy in (\ref{PV})  is  correct (see, for example
\citep{BHJ92}) and is equal to
$C\sum_1^np_i^2(1\vee\sum_1^np_i)^{-1}$. Similarly, in (\ref{PVA})
the order of accuracy is $C(\max p_i)^2$. As one can expect, the
accuracy of approximation is trivial, if all $p_i$ are uniformly
bounded from zero, i.e., $p_i>C$.
The accuracy of approximation is much better for $\G$.
\begin{theorem} \label{G2T} Let conditions (\ref{nu12}) and (\ref{lambda}) be
satisfied. Then, for all $n$,
\begin{eqnarray}
\norm{F_n-\G}&\leqslant&
C R_1\big\{1+\Gamma_1\min(1,\lambda^{-1})\big\}\min(1,\lambda^{-3/2}),\label{GV}
\\
\norm{F_n-\G(\dirac+\Gamma_3(\dirac_1-\dirac)^3)}&\leqslant&
C\big\{1+\Gamma_1\min(1,\lambda^{-1})\big\}\big(R_1^2\min(1,\lambda^{-3})\nonumber\\&&+R_2\min(1,\lambda^{-2})\big)\label{GVA}.
\end{eqnarray}
\end{theorem}
If, instead of (\ref{lambda}), we assume (\ref{3ab}), then $\lambda\geqslant
C\Gamma_1$  and
$1+\Gamma_1\min(1,\lambda^{-1})\leqslant C$. If, in addition, all $X_i$ do not depend on $n$ and
are bounded, then estimates in (\ref{GV}) and (\ref{GVA}) are of
 orders $O(n^{-1/2})$ and $O(n^{-1})$, respectively. Thus, the
order of accuracy is comparable to CLT and Edgeworth's expansion.
 If all
$X_i\sim Be(1,p_i)$ are independent, then the order of accuracy in
(\ref{GV})  is the right one (see \citep{Kru86}) and is equal to
$C\sum_1^np_i^3(1\vee\sum_1^np_i)^{-3/2}$.

Approximation $\G$ has two parameters, but: a) is not always a
distribution, b) its "probabilities" are not easily calculable.
Some authors argue (see, for example, \citep{BrX01}) that,
therefore, probabilistic approximations are more preferable. We start from the negative binomial approximation. Observe, that
the negative binomial approximation is meaningful only if $\Var
S_n
> \Expect S_n$.


\begin{theorem} \label{NBtheorem} Let conditions (\ref{nu12}) and (\ref{3ab}) be
satisfied and let $\Gamma_2>0$. Then, for all $n$,
\begin{eqnarray}
\norm{F_n-\NB(r,\bar{q})}&\leqslant&C \min(1,\Gamma_1^{-3/2})(R_1+\Gamma_2^2\Gamma_1^{-1}),\hskip 2cm\label{NBV}\\
\norm{F_n-\NB(r,\bar{q})\big(\dirac+[\Gamma_3-4\Gamma_2^2(3\Gamma_1)^{-1}](\dirac_1-\dirac)^3\big)}
&\leqslant&C\big\{
R_1^2\min(1,\Gamma_1^{-3})\nonumber\\
\lefteqn{+R_2\min(1,\Gamma_1^{-2})+\Gamma_2^2\Gamma_1^{-1}\ab{\Gamma_3-4\Gamma_2^2(3\Gamma_1)^{-1}}\min(1,\Gamma_1^{-3})
+\Gamma_2^{3}\Gamma_1^{-2}\min(1,\Gamma_1^{-2})\big\}.}\hskip 6cm
\label{NBVA}
\end{eqnarray}
\end{theorem}
It seems that asymptotic expansion for the negative binomial
approximation was so far never considered in the context of
1-dependent summands. If all $X_i$ do not depend on $n$ and are
bounded, the accuracies of approximation in (\ref{NBV}) and
(\ref{NBVA}) are $O(n^{-1/2})$ and $O(n^{-1})$, respectively.

If  $\Var S_n < \Expect S_n$, it is more natural to use the
binomial approximation.

\begin{theorem} \label{Bintheorem} Let conditions (\ref{nu12}) and (\ref{3ab}) be
satisfied,  $\Gamma_1\geqslant 1$ and $\Gamma_2<0$. Then, for all
$n$,
\begin{eqnarray}
\norm{F_n-\BI(N,\bar{p})}&\leqslant&C (\Gamma_2^2\Gamma_1^{-5/2}+
R_1\Gamma_1^{-3/2}),\label{BIV}\\
\norm{F_n-\BI(N,\bar{p})\big(\dirac+
[\Gamma_3-N\bar{p}^3/3](\dirac_1-\dirac)^3\big)}
&\leqslant&C\big\{ R_1^2\Gamma_1^{-3}+R_2\Gamma_1^{-2}
+\ab{\Gamma_2}^{3}\Gamma_1^{-4}\nonumber\\
&&+\Gamma_2^2\Gamma_1^{-3}
+\Gamma_2^2\Gamma_3\Gamma_1^{-4}\big\}. \label{BIVA}
\end{eqnarray}
\end{theorem}
If all the $X_i$ do not depend on $n$ and are bounded, the accuracies
of approximation in (\ref{BIV}) and (\ref{BIVA})  are
$O(n^{-1/2})$ and $O(n^{-1})$, respectively.

 In this paper we consider the total variation norm only. It must be noted that  formula of inversion for probabilities allows to prove local estimates too. If $\lambda>1$, then   local estimates are equal to (\ref{PV}) -- (\ref{BIVA}) multiplied by factor $\lambda^{-1/2}$.

\section{Applications}

\textbf{1. Asymptotically sharp constant for the negative binomial
approximation to 2-runs.} As already mentioned  above, the 2-
runs statistic is one of the best investigated cases of sums of
1-dependent discrete random variables. It is easy to check that
the rate of accuracy in (\ref{BXNB}) is $O(pn^{-1/2})$. However,
the constant 64.4 is not particularly small. Here, we shall
show, that, on the other hand, \emph{asymptotically} sharp
constant is small. Asymptotically sharp constant  can be used heuristically to get the  idea about the
 magnitude of  constant in (\ref{NBV}). We shall consider 2-runs with edge effects,
which we think to be more realistic case than $S^{*}$. Let
$S_\xi=\xi_1+\xi_2+\dots+\xi_n$, where $\xi_i=\eta_i\eta_{i+1}$
and $\eta_i\sim Be(p)$, $(i=1,2,\dots,n+1)$ are independent
Bernoulli variables. The sum $S^{*}$ differs from $S_\xi$  by the
last summand only, which is equal to $\eta_n\eta_1$. As shown in \citep{PeCe1}, for $S_\xi$
we have
\[
\Gamma_1=np^2,\quad\Gamma_2= \frac{np^3(2-3p)-2p^3(1-p)}{2},\quad
\Gamma_3=\frac{np^4(3-12p+10p^2)-6p^4(1-p)(1-2p)}{3}.
\]
 Let $\NB(r,\bar{q})$ be defined as in
(\ref{nbinom}) and
\[
\tilde
C_{TV}=\frac{1}{3}\sqrt{\frac{2}{\pi}}(1+4\ee^{-3/2})=0.5033...~
\]

\begin{theorem} \label{NBASHARP} Let $p\leqslant 1/20$, $np^2\geqslant 1$. Then
\[\Ab{ \norm{\eL(S_\xi)-\NB(r,\bar{q})}-\tilde
C_{TV}\frac{p}{\sqrt{n}}}\leqslant C\bigg(\frac{p^2}{\sqrt{n}}+\frac{1}{n}\bigg).\label{ASNBT}
\]
\end{theorem}

\noindent We now get the following corollary.
\begin{corollary} \label{corasharp}
 Let $p\to 0$ and $np^2\to \infty$, as $n \to \infty$. Then
 \[
\lim_{n\to\infty}\frac{\norm{\eL(S_\xi)-\NB(r,\bar{q})}\sqrt{n}}{p}=\tilde
C_{TV}.
\]
\end{corollary}
\vskip 0.5cm

 \textbf{2. Binomial approximation to $N(k_1,k_2)$ events.} Let $\eta_i\sim
 Be(p)$,($0<p<1$)
be independent Bernoulli variables and let
$Y_j=(1-\eta_{j-m+1})\cdots
(1-\eta_{j-k_2})\eta_{j-k_2+1}\cdots\eta_{j-1}\eta_j$,
$j=m,m+1,\dots,n$, $k_1+k_2=m$. Further, we assume that $k_1>0$
and $k_2>0$. Let $N(n;k_1,k_2)=Y_m+Y_{m+1}+\cdots+Y_n$. We denote
the distribution of $N(n;k_1,k_2)$ by $\Ha$. Let
$a(p)=(1-p)^{k_1}p^{k_2}$. It is well known that $N(n;k_1,k_2)$
has limiting Poisson distribution and  the accuracy of Poisson
approximation is $O(a(p))$, see
 \citep{HuTs91} and \citep{Vel04}, respectively. However, Poisson approximation
has just one parameter. Consequently, the closeness of $p$ to zero
is crucial. We can expect any two-parametric approximation to be
more universal. It is proved in  \citep{Up09} that
 \begin{eqnarray*} \Expect
N(n;k_1,k_2)&=&(n-m+1)a(p),\\
\Var N(n;k_1,k_2)&=&(n-m+1)a(p)+(1-4m+3m^2-n(2m-1))a^2(p).
\end{eqnarray*}
 Under quite mild assumptions  $\Var
N(n;k_1,k_2)<\Expect N(n;k_1,k_2)$. Consequently, the natural
probabilistic approximation is binomial one. The binomial
approximation to $N(n;k_1,k_1)$ was already considered  in
\citep{Up09}. Regrettably, the estimate in \citep{Up09} contains
expression which is of the constant order when $a(p)\to 0$.

Note that $Y_1,Y_2,\dots$ are $m$-dependent. Consequently, results
of the previous Section can not be applied directly. However, one
can group summands in the following natural way:
\[N(n;k_1,k_2)=(Y_m+Y_{m+1}+\cdots+Y_{2m-1})+(Y_{2m}+Y_{2m+1}+\cdots+Y_{3m-1})+\dots=X_1+X_2+\dots\]
Each $X_j$, with probable exception of the last one, contains $m$
 summands. It is not difficult  to check that $X_1,X_2,\dots$ are
 1-dependent
 Bernoulli variables.
All parameters can be written explicitly. Set $N=\lfloor \tilde
N\rfloor$ be the integer part of $\tilde N$,
\[ \tilde
N=\frac{(n-m+1)^2}{(n-m+1)(2m-1)-m(m-1)},\quad \tilde N=
N+\epsilon, \quad
0\leqslant\epsilon<1,\quad\bar{p}=\frac{(n-m+1)a(p)}{N}.\] For the
asymptotic expansion we need the following notation
\begin{equation*}
A=\frac{a^3(p)}{6}(n-m+1)m(m-1).
\end{equation*}

The two-parametric binomial approximation is more natural, when
$\Expect N(n;k_1,k_2)\geqslant 1$, which means that we deal with
large values of $n$ only.

\begin{theorem} \label{k1k2th} Let $(n-m+1)a(p)\geqslant 1$ and $ma(p)\leqslant 0.01$. Then
\begin{eqnarray}
\norm{\Ha-\BI(N,\bar{p})}&\leqslant&C
\frac{a^{3/2}(p)m^2}{\sqrt{n-m+1}},\label{BIVk1k2}\\
\norm{\Ha-\BI(N,\bar{p})\big(\dirac+ A(\dirac_1-\dirac)^3\big)}
&\leqslant&C\frac{a(p)m^2(a(p)m+1)}{n-m+1}. \label{BIVAk1k2}
\end{eqnarray}
\end{theorem}

Note that the assumption $ma(p)\leqslant 0.01$  in Theorem
\ref{k1k2th} is not very restrictive on $p$ when $k_1,k_2
> 1$. For example, it is satisfied for $p\leqslant 1/4$ and
$N(n;4,4)$.
\begin{theorem} \label{BINSHARP} Let $(n-m+1)a(p)\geqslant 1$ and $ma(p)\leqslant 0.01$. Then
\begin{eqnarray*}
\lefteqn{\Ab{ \norm{\Ha-\BI(N,\bar{p})}- \tilde
C_{TV}\frac{a^{3/2}(p)m(m-1)}{2\sqrt{n-m+1}}}}\hskip 4cm\nonumber\\
&\leqslant& C(m)\frac{a^{3/2}(p)m(m-1)}{\sqrt{n-m+1}}
\bigg(\frac{1}{\sqrt{(n-m+1)a(p)}}+\frac{1}{\tilde
N-1}+a(p)\bigg).\label{ASTVk1k2}
\end{eqnarray*}
\end{theorem}
Constant $C(m)$ depends on $m$.
\begin{corollary} \label{corasharp}
 Let $m$ be fixed, $a(p)\to 0$, $(n-m+1)a(p)\to \infty$, as $n\to\infty$. Then
 \begin{eqnarray*}
\lim_{n\to\infty}\frac{\norm{\Ha-\BI(N,\bar{p})}\sqrt{n-m+1}}{a^{3/2}(p)m(m-1)}&=&\frac{\tilde
C_{TV}}{2}.
\end{eqnarray*}
\end{corollary}


\section{Auxiliary results}
In this section, some auxiliary results from other papers are collected.
For the sake of brevity, further we will use the notation $U=\dirac_1-\dirac$.
First, we need representation of the characteristic function $\w F(t)$ as
product of functions.

\begin{lemma}\label{Hversion} Let conditions (\ref{nu12}) and (\ref{lambda}) be
satisfied. Then
\[\w F(t)=\vfi_1(t)\vfi_2(t)\dots\vfi_n(t),\]
 where $\vfi_1(t)=\Expect e^{\ii t X_1}$ and, for $k=2,\dots, n$,
\[\vfi_k(t)=1+ \Expect (e^{\ii t X_k}-1)+\sum_{j=1}^{k-1}\frac{\wE((e^{\ii
tX_j}-1), (e^{\ii t X_{j+1}}-1),\dots (e^{\ii t
X_k}-1))}{\vfi_j(t)\vfi_{j+1}(t)\dots \vfi_{k-1}(t)}.\]
\end{lemma}

Lemma \ref{Hversion} follows from more general Lemma
3.1 in \citep{H82}. Representation holds for all $t$, since the
assumption of Lemma 3.1
\[
\sqrt{\Expect\ab{\ee^{\ii t X_k}-1}^2}\leqslant \sqrt{2
\Expect\ab{\ee^{\ii t X_k}-1}}\leqslant \sqrt{2\nu_1(k)}\leqslant
\sqrt{0.02}<1/6
\]
is satisfied for all $t$.

\begin{lemma} \label{a10}
Let $t\in(0,\infty)$, $0<p<1$ and $n,j=1,2,\dots$. We then  have
\begin{eqnarray*}
\norm{U^2\ee^{tU}}&\leqslant&\frac{3}{t\ee},\qquad
\norm{U^j\ee^{tU}}\leqslant\Big(\frac{2j}{t\ee}\Big)^{j/2},
\quad
\norm{U^j(\dirac+pU)^n}
\leqslant\binomial{n+j}{j}^{-1/2}(p(1-p))^{-j/2}.
\end{eqnarray*}
\end{lemma}
The first inequality was proved in \citep{Roo01} (formula~(29)).
The second bound follows from formula (3.8) in \citep{DP88} and
the properties of the total variation norm.
 For the proof of the third estimate,  see Lemma~4 from \citep{Roo00}.

%

\begin{lemma}\label{sharpC}
Let $t>0$ and $p\in(0,1)$. Then
\begin{equation*}
\Ab{\norm{U^3\ee^{tU}}-\frac{3 \tilde C_{TV}}{t^{3/2}}}
\leqslant\frac{C}{t^{2}},\qquad
\Ab{\norm{U^3(\dirac+pU)^n} -\frac{3\tilde
C_{TV}}{(np(1-p))^{3/2}} }\leqslant\frac{C}{(np(1-p))^{2}}.
\end{equation*}
\end{lemma}
The statements in Lemma \ref{sharpC} follow from a more general Proposition 4 in
\citep{Roo99} and from \citep{CeRo06}.
\begin{lemma} \label{simple} Let $\lambda>0$ and $k=0,1,2,\dots$.
Then
\[
\ab{\sin(t/2)}^k\ee^{-\lambda\sin^2(t/2)}\leqslant
\frac{C(k)}{\lambda^{k/2}},\quad
\int_{-\pi}^\pi \ab{\sin(t/2)}^k\ee^{-\lambda\sin^2(t/2)}\dd
t\leqslant \frac{C(k)}{\max(1,\lambda^{(k+1)/2})}.
\]
\end{lemma}
Both estimates are  trivial. Note that, for
$\ab{t}\leqslant\pi$, we have $\ab{\sin(t/2)}\geqslant
\ab{t}/\pi$.


\begin{lemma} \label{apvertimas} Let $M$ be   finite variation measure concentrated on
integers, $\sum_k\ab{k}\ab{M\{k\}}<\infty$. Then for any  $v\in\RR$ and $u
>0$ the following inequality is valid
\begin{equation}
 \norm{M}\leqslant\bigl(1+u\pi\bigr)^{1/2}
\biggl(\frac{1}{2\pi}\int\limits_{-\pi}^{\pi}\ab{\w M(t)}^2+
\frac{1}{u^2}\ab{\bigl(\hbox{\rm e}^{-\ii tv}\w M(t)\bigr)'}^2
\hbox{\rm d}t\biggr)^{1/2}. \label{TVAP}
 \end{equation}
\end{lemma}

The estimate (\ref{TVAP}) is well-known; see, for example,
\citep{Pre85}.


\begin{lemma} (\citep{Berg51}) For all numbers $A,B>0$, $s=0,1,2,\dots,n$, the
following identity holds:
\begin{equation}
A^n\,=\,\sum_{m=0}^s\binomial{n}{m}B^{n-m}(A-B)^m+\sum_{m=s+1}^n\binomial{m-1}{s}A^{n-m}(A-B)^{s+1}B^{m-s-1}.
            \label{bergident}
\end{equation}
\end{lemma}



\begin{lemma} \label{factorial} Let $s=1,2,3$. For all $t\in\RR$,
\begin{eqnarray*}
\Expect \exponent{\ii t X_k} &=&
1+\sum_{l=1}^s\nu_l(k)\frac{(\eit-1)^l}{l!}+\theta\nu_{s+1}(k)\frac{\ab{\eit-1}^{s+1}}{s!},\\
\Expect (\exponent{\ii t X_k})' &=&
\sum_{l=1}^s\nu_l(k)\frac{\ii\eit(\eit-1)^{l-1}}{l!}+\theta\nu_{s+1}(k)\frac{\ab{\eit-1}^{s}}{(s-1)!}.
\end{eqnarray*}
\end{lemma}
 Lemma \ref{factorial} is a particular case of Lemma 3 from
\citep{SC88}.

\begin{lemma}(\citep{H82}) \label{Hei3aa} Let $Z_1,Z_2,\dots,Z_k$ be 1-dependent complex-valued random variables with
$\Expect\ab{Z_m}^2<\infty$,  $1 \leqslant m \leqslant k. $ Then
\[
\ab{\wE (Z_1, Z_2, \cdots, Z_k)}\leqslant
2^{k-1}\prod_{m=1}^k(\Expect\ab{Z_m}^2)^{1/2}.
\]
\end{lemma}

\section{Preliminary results}

 Let
$z=\eit-1$ and $Z_j=\exponent{\ii t X_j}-1$. As before we assume
that $\nu_j(k)=0$ and $X_k=0$ for $k\leqslant 0$.
 Also, we omit the argument $t$, wherever possible and,  for example,
write $\varphi_k$ instead of $\varphi_k(t)$.


The next lemma can easily be proved by induction.
\begin{lemma}\label{Lema1} For all $t\in\RR$ and $k\geqslant 2$, the following estimate holds:
\begin{equation}
\wE^{+}(\ab{Z_1},\dots, \ab{Z_k})\leqslant
4\wE^{+}(\ab{Z_1},\dots,\ab{Z_{k-1}}). \label{shorting}
\end{equation}
\end{lemma}

\begin{lemma} \label{shortphi} Let $\max_k\nu_1(k)\leqslant 0.01$.
Then, for $k=1,2,\dots,n$,
\begin{eqnarray}
\ab{\vfi_k-1}&\leqslant&\frac{1}{10},\quad
\frac{1}{\ab{\vfi_k}}\leqslant\frac{10}{9},\label{sh1}\\
\quad\ab{\vfi_k-1}&\leqslant& \ab{z}[(0.66)\nu_1(k-1)+(4.13)\nu_1(k)],
\label{sh2}\\
\ab{\vfi_k-1-\Expect Z_k}&\leqslant& \sin^2(t/2)[
(0.374)\nu_1(k)+(0.288)\nu_1(k-1)\nonumber\\&&+(15.58)\Expect X_{k-1}X_k
+(0.1)\Expect X_{k-2}X_{k-1}].\label{sh3}
\end{eqnarray}
\end{lemma}

\noindent \Proof  We repeatedly apply below the following trivial
inequalities:
\begin{equation}
\ab{z}\leqslant 2,\qquad\ab{Z_k}\leqslant 2,\qquad
\ab{Z_k}\leqslant X_k\ab{z}. \label{trivial}
\end{equation}
The second estimate in (\ref{sh1}) follows  from the first
estimate:
\[ \ab{\vfi_k}\geqslant\ab{1-\ab{\vfi_k-1}}\geqslant
1-(1/10)=9/10.\]
 The first estimate in  (\ref{sh1}) follows from (\ref{sh2})
 and  (\ref{trivial}) and by the assumption of  the lemma.
 It remains to prove the (\ref{sh2}) and (\ref{sh3}). Both proofs are very similar.  From Lemma \ref{Hversion} and
equation (\ref{sh1}), we get
\begin{eqnarray}
\ab{\vfi_k-1-\Expect
Z_k}&\leqslant&\frac{\ab{\wE(Z_{k-1},Z_k)}}{\ab{\vfi_{k-1}}}+
\frac{\ab{\wE(Z_{k-2},Z_{k-1},Z_k)}}{\ab{\vfi_{k-3}\vfi_{k-1}}}
+\frac{\ab{\wE(Z_{k-3},\dots,Z_k)}}{\ab{\vfi_{k-3}\cdots\vfi_{k-1}}}\nonumber\\
&&+\frac{\ab{\wE(Z_{k-4},\dots,Z_k)}}{\ab{\vfi_{k-4}\cdots\vfi_{k-1}}}+
\frac{\ab{\wE(Z_{k-5},\dots,Z_k)}}{\ab{\vfi_{k-5}\cdots\vfi_{k-1}}}+
 \sum_{j=1}^{k-6}\frac{\ab{\wE(Z_j,\dots,Z_k)}}{\ab{\vfi_j\vfi_{j+1}\cdots\vfi_{k-1}}}\nonumber
 \end{eqnarray}
 \begin{eqnarray}
&\leqslant&\bigg(\frac{10}{9}\bigg)\ab{\wE(Z_{k-1},Z_k)}+
\bigg(\frac{10}{9}\bigg)^2\ab{\wE(Z_{k-2},Z_{k-1},Z_k)}
+\bigg(\frac{10}{9}\bigg)^3\ab{\wE(Z_{k-3},\dots,Z_k)}\nonumber\\
&&+\bigg(\frac{10}{9}\bigg)^4\ab{\wE(Z_{k-4},\dots,Z_k)}+
\bigg(\frac{10}{9}\bigg)^5\ab{\wE(Z_{k-5},\dots,Z_k)}\nonumber\\
&&+
 \sum_{j=1}^{k-6}\bigg(\frac{10}{9}\bigg)^{k-j}\ab{\wE(Z_j,\dots,Z_k)}.
\label{sh4}
\end{eqnarray}
By (\ref{trivial}) and  Lemma \ref{Hei3aa}, we obtain
\begin{equation}
\Expect\ab{Z_j}\leqslant \nu_1(j)\ab{z}\leqslant
0.02\ab{\sin(t/2)},\quad \Expect\ab{Z_j}^2\leqslant
2\Expect\ab{Z_j}\leqslant 2\nu_1(j)\ab{z}=4\nu_1(j)\ab{\sin(t/2)}
\label{ez1} \end{equation} and
\begin{eqnarray}
\ab{\wE(Z_j,\dots Z_k)}&\leqslant&
2^{k-j}2^{(k-j+1)/2}\ab{z}^{(k-j+1)/2}\prod_{l=j}^k\sqrt{\nu_1(l)}\nonumber\\
&\leqslant&
2^{2(k-j)-1}\ab{z}^2\sqrt{\nu_1(k)\nu_1(k-1)}0.1^{k-j-1}\nonumber\\
&\leqslant&
4^{k-j}\sin^2\frac{t}{2}[\nu_1(k)+\nu_1(k-1)]0.1^{k-j-1}\nonumber\\
&&=
10\sin^2\frac{t}{2}[\nu_1(k)+\nu_1(k-1)](0.4)^{k-j}.\label{ez00}
\end{eqnarray}
Consequently,
 \begin{eqnarray}
 \sum_{j=1}^{k-6}\bigg(\frac{10}{9}\bigg)^{k-j}\ab{\wE(Z_j,\dots,Z_k)}
&\leqslant&
10\sin^2\frac{t}{2}[\nu_1(k)+\nu_1(k-1)]\sum_{j=1}^{k-6}\bigg(\frac{4}{9}\bigg)^{k-j}\nonumber\\
&\leqslant& 2\sin^2\frac{t}{2}[\nu_1(k)+\nu_1(k-1)](0.0694).
\label{sh5}
\end{eqnarray}
By 1-dependence,  (\ref{trivial}) and H\"older's inequality (see
also \citep{H82}), we have for $j\geqslant 3$,
\begin{eqnarray}
\ab{\Expect Z_{k-j}\cdots
Z_k}&\leqslant&\prod_{i=k-j}^k\sqrt{\Expect\ab{Z_i}^2}
\leqslant\prod_{i=k-j}^k\sqrt{2\nu_1(i)\ab{z}} \leqslant
2^{(j+1)/2}\ab{z}^{(j+1)/2}\sqrt{\nu_1(k-1)\nu_1(k)}(0.1)^{j-1}\nonumber\\
&\leqslant&
2^{j-1}\ab{z}^2\frac{\nu_1(k-1)+\nu_1(k)}{2}(0.1)^{j-1}\nonumber\\
&=& 2^j \sin^2(t/2)[\nu_1(k-1)+\nu_1(k)](0.1)^{j-1}. \label{ezkj}
\end{eqnarray}
Moreover, for any $j$,
\begin{equation}
\ab{\Expect Z_{j-1}Z_j}\leqslant 2\Expect\ab{Z_j}\leqslant
4\ab{\sin(t/2)}\nu_1(j),\quad \ab{\Expect Z_{j-1}Z_j}\leqslant
\ab{z}^2\Expect{X_{j-1}X_j}=4\sin^2(t/2)\Expect{X_{j-1}X_j}\label{ez2}
\end{equation}
and
\begin{equation}
\ab{\Expect Z_{j-2}Z_{j-1}Z_j}\leqslant
2\Expect\ab{Z_{j-1}Z_j}\leqslant 8\sin^2(t/2)\Expect X_{j-1}X_j.
\label{ez3}
\end{equation}
Therefore, from  (\ref{ez1}), we have
\begin{equation}
\ab{\wE(Z_{j-1},Z_j)}\leqslant\Expect\ab{Z_{j-1}Z_j}+\nu_1(j-1)\nu_1(j)\ab{z}^2\leqslant
2.02\ab{z}\nu_1(j)\leqslant 0.0404\ab{\sin(t/2)}. \label{rough1}
\end{equation}
Similarly, applying (\ref{ezkj}), (\ref{ez2}),  (\ref{ez3}) and
(\ref{wez2}), we obtain the following rough estimates:
\begin{eqnarray}
\ab{\wE(Z_{j-2},Z_{j-1},Z_j)}&\leqslant&
\ab{z}\ab{\sin(t/2)}\{0.2\nu_1(j-1)+0.2804\nu_1(j)\}\leqslant
0.01\sin^2(t/2),\nonumber\\
\ab{\wE(Z_{j-3},\dots,Z_j)}&\leqslant&
\ab{z}\ab{\sin(t/2)}\{0.044\nu_1(j-1)+0.1348\nu_1(j)\}\leqslant
0.0036\sin^2(t/2),\label{rough3}\\
\ab{\wE(Z_{j-4},\dots,Z_j)}&\leqslant&
\ab{z}\ab{\sin(t/2)}\{0.0169\nu_1(j-1)+0.0405\nu_1(j)\}\leqslant
0.00115\sin^2(t/2). \nonumber
\end{eqnarray}
Taking into account that $\nu_1(k-1)\leqslant 0.01$, we get
\begin{equation}
\ab{\wE(Z_{k-1},Z_k)}\leqslant\Expect\ab{
Z_{k-1}Z_k}+\nu_1(k-1)\nu_1(k)\ab{z}^2\leqslant
\sin^2(t/2)\{4\Expect X_{k-1}X_k+0.04\nu_1(k)\}.\label{wez2}
\end{equation}
Similarly, taking into account (\ref{ezkj})--(\ref{rough3}), we
get
\begin{eqnarray}
\ab{\wE(Z_{k-2},Z_{k-1},Z_k)}&\leqslant&\sin^2(t/2)\{8.08\Expect
X_{k-1}X_k+0.08\Expect X_{k-2}X_{k-1}+0.0008\nu_1(k-1)\},
\nonumber\\
\ab{\wE(Z_{k-3},\dots,Z_k)}&\leqslant&\sin^2(t/2)\{0.3216\Expect
X_{k-1}X_k+0.08\nu_1(k-1)+0.1\nu_1(k)\},
\nonumber\\
\ab{\wE(Z_{k-4},\dots,Z_k)}&\leqslant&\sin^2(t/2)\{0.3632\Expect
X_{k-1}X_k+0.0176\nu_1(k-1)+0.0248\nu_1(k)\},
\nonumber\\
\ab{\wE(Z_{k-5},\dots,Z_k)}&\leqslant&\sin^2(t/2)\{0.0944\Expect
X_{k-1}X_k+0.0068\nu_1(k-1)+0.0091\nu_1(k)\}. \label{wez6}
\end{eqnarray}
Combining (\ref{sh5}), (\ref{wez2})--(\ref{wez6}) with (\ref{sh4})
we prove (\ref{sh3}).

For the proof of (\ref{sh2}), we
apply
 mathematical induction. Let us assume that (\ref{sh1}) holds for first $k-1$
functions and let $k\geqslant 6$. Then the proof is almost
identical to the proof of (\ref{sh3}). We expand $\vfi_k$ just
like in (\ref{sh4}):
\[
\ab{\vfi_k-1} \leqslant \Expect\ab{ Z_k}+
\bigg(\frac{10}{9}\bigg)\ab{\wE(Z_{k-1},Z_k)}+\dots+
 \sum_{j=1}^{k-4}\bigg(\frac{10}{9}\bigg)^{k-j}\ab{\wE(Z_j,\dots,Z_k)}.
\]
Applying (\ref{ezkj}), (\ref{ez1}) and  (\ref{rough1})--(\ref{rough3}), we easily complete the proof of (\ref{sh2}).
The proof
for $k<6$ is analogous. \qed

\begin{lemma}\label{lambdak} Let  $\nu_1(k)\leqslant 0.01$, $\nu_2(k)<\infty$, for $ 1 \leq k \leq n$.
Then, for all $t\in\RR$,
\[
\ab{\vfi_k}\leqslant
1-\lambda_k\sin^2(t/2)\leqslant\exponent{-\lambda_k\sin^2(t/2)}\]
and
\[\ab{\w F(t)}\leqslant\prod_{k=1}^n\ab{\vfi_k}\leqslant
\exponent{-1.3\lambda\sin^2(t/2)}.\]
 Here $\lambda_k= 1.606\nu_1(k)-0.288\nu_1(k-1)-2\nu_2(k)-0.1\Expect
 X_{k-2}X_{k-1}-15.58\Expect X_{k-1}X_{k}$ and $\lambda$ is
 defined by (\ref{lambda}).
 \end{lemma}

\noindent \Proof  From Lemma \ref{factorial}  it follows that
\[1+\Expect Z_k=\Expect\exponent{\ii t X_k}=1+\nu_1(k)z+\theta\frac{\nu_2(k)\ab{z}^2}{2}.\]
Therefore
\begin{equation*}
\ab{\vfi_k}\leqslant\ab{1+\Expect Z_k}+\ab{\vfi_k-1-\Expect
Z_k}\leqslant
\ab{1+\nu_1(k)z}+\frac{\nu_2(k)}{2}\ab{z}^2+\ab{\vfi_k-1-\Expect
Z_k}. \label{l1}
\end{equation*}
Applying the definition of the square of the absolute value for
complex number we get
\[\ab{1+\nu_1(k)z}^2=(1-\nu_1(k)\cos t)^2+(\nu_1(k)\sin
t)^2=1-4\nu_1(k)(1-\nu_1(k))\sin^2(t/2).\]
 Consequently,
 \[\ab{1+\nu_1(k)z}\leqslant\sqrt{ 1-4\nu_1(k)(1-\nu_1(k))\sin^2(t/2)}\leqslant
 1-2\nu_1(k)(1-\nu_1(k))\sin^2(t/2).\]
Combining the last estimate with (\ref{sh3}), we get the first estimate of the lemma. The second
estimate follows immediately. \qed

For expansions of $\vfi_k$ in powers of $z$, we use the following
notation:
\begin{eqnarray*}
\gamma_2(k)&=&\frac{\nu_2(k)}{2}+\wE(X_{k-1},X_k),\nonumber\\
\gamma_3(k)&=&\frac{\nu_3(k)}{6}+\frac{\wE_2(X_{k-1},X_k)}{2}+\wE(X_{k-2},X_{k-1},X_k)-\nu_1(k-1)\wE(X_{k-1},X_k),\nonumber\\
r_0(k)&=&\nu_2(k)+\sum_{l=0}^3\nu_1^2(k-l)+\Expect
X_{k-1}X_k,\label{r0k}\\
r_1(k)&=&\nu_3(k)+\sum_{l=0}^5\nu_1^3(k-l)+\nu_1(k-1)\Expect
X_{k-1}X_k+\wE_2^{+}(X_{k-1},X_k)\nonumber\\
&&+\wE^{+}(X_{k-2},X_{k-1},X_k),\label{r1k}
\\
r_2(k)&=&\nu_4(k)+\sum_{l=0}^7\nu_1^4(k-l)+\nu_2^2(k)+\nu_2^2(k-1)+(\Expect
X_{k-1}X_k)^2+ (\Expect
X_{k-2}X_{k-1})^2\nonumber\\
&&+(\nu_1(k-2)+\nu_1(k-1))\wE_2^{+}(X_{k-1},X_k)+\sum_{l=0}^3\nu_1(k-l)\wE^{+}(X_{k-2},X_{k-1},X_k)\nonumber\\
&&+\wE_2^{+}(X_{k-2},X_{k-1},X_k)+\wE_3^{+}(X_{k-1},X_{k})+\wE^{+}(X_{k-3},X_{k-2},X_{k-1},X_k).\label{r2k}
\end{eqnarray*}

\begin{lemma}\label{f} Let condition (\ref{nu12}) be satisfied, $k=1,\dots,n$.
Then, for all $t\in\RR$,
\begin{eqnarray}
\vfi_k&=&1+\nu_1(k)z+\theta C\ab{z}^2r_0(k),\label{f1}\\
\vfi_k&=&1+\nu_1(k)z+\gamma_2(k)z^2+\theta C\ab{z}^3r_1(k),\label{f2}\\
\vfi_k&=&1+\nu_1(k)z+\gamma_2(k)z^2+\gamma_3(k)z^3+\theta C\ab{z}^4r_2(k),\label{f3}
\end{eqnarray}\begin{eqnarray}
\frac{1}{\vfi_{k-1}}&=&1+C\theta\ab{z}\{\nu_1(k-2)+\nu_1(k-1)\},\label{f4}\\
\frac{1}{\vfi_{k-1}}&=&1-\nu_1(k-1)z+C\theta\ab{z}^2\big\{\nu_2(k-1)+\sum_{l=1}^4\nu_1^2(k-l)+\Expect
X_{k-2}X_{k-1}\big\},\label{f5} \\
\frac{1}{\vfi_{k-1}}&=&1-\nu_1(k-1)z-\bigg(\frac{\nu_2(k-1)}{2}-\nu_1^2(k-1)+\wE(X_{k-2},X_{k-1})\bigg)z^2\nonumber\\
&&+C\theta\ab{z}^3\Big\{\nu_3(k-1)+\sum_{l=1}^6\nu_1^3(k-l)+\wE_2^{+}(X_{k-2},X_{k-1})\nonumber\\
&&+\wE^{+}(X_{k-3},X_{k-2},X_{k-1})\nonumber\\
&&+(\nu_1(k-2)+\nu_1(k-1))[\nu_2(k-1)+\Expect
X_{k-2}X_{k-1}]\Big\},
\label{f6}\\
\frac{1}{\vfi_{k-1}\vfi_{k-2}}&=&1+C\theta\ab{z}\{\nu_1(k-3)+\nu_1(k-2)+\nu_1(k-1)\},\label{f7}\\
(\vfi_k-1)^2&=&\nu_1^2(k)z^2+C\theta\ab{z}^3[\nu_1(k-1)\nonumber\\
&&+\nu_1(k)]\Big\{\nu_2(k)+\sum_{l=0}^3\nu_1^2(k-l)+\Expect
X_{k-1}X_k\Big\},
\label{f8}\\
(\vfi_k-1)^3&=&\nu_1^3(k)z^3\nonumber\\
&&+C\theta\ab{z}^4\Big\{\nu_2^2(k)+\sum_{l=0}^3\nu_1^4(k-l)+(\Expect
X_{k-1}X_k)^2\Big\}.\label{f9}
\end{eqnarray}
\end{lemma}

\noindent \Proof Further on we assume that $k\geqslant 7$. For smaller
values of $k$, all proofs just become shorter. The lemma is proved in
four steps. First, we prove (\ref{f1}), (\ref{f2}), (\ref{f4}) and
(\ref{f7}). Second, we obtain (\ref{f5}) and (\ref{f8}). Then we
prove (\ref{f6}) and (\ref{f9}). The final step is the proof of
(\ref{f3}). At each step, we employ results from the previous step.
Since all proofs are very similar, we give just some of them.

\noindent Due to (\ref{sh1}), we have
\[\frac{1}{\vfi_{k-1}}=\frac{1}{1-(1-\vfi_{k-1})}=1+\sum_{j=1}^\infty(1-\vfi_{k-1})^j=1
+C\theta\ab{1-\vfi_{k-1}}.
\]
Therefore, (\ref{f4}) and (\ref{f7}) follow from (\ref{sh2}) and (\ref{sh1}).\\
 From Lemmas \ref{Hversion}, \ref{factorial}, \ref{Lema1},
equation (\ref{sh1}) and second estimate in (\ref{ez00}), we get
\begin{eqnarray*}
\ab{\vfi_k}&=&1+\Expect
Z_k+\frac{\ab{\wE(Z_{k-1},Z_k)}}{\ab{\vfi_{k-1}}}+
\frac{\ab{\wE(Z_{k-2},Z_{k-1},Z_k)}}{\ab{\vfi_{k-2}\vfi_{k-1}}}
 +\sum_{j=1}^{k-3}\frac{\ab{\wE(Z_j,\dots,Z_k)}}{\ab{\vfi_j\vfi_{j+1}\cdots\vfi_{k-1}}}\nonumber\\
&\leqslant&
1+\nu_1(k)z+C\theta\ab{z}^2\nu_2(k)+C\theta\wE^{+}(\ab{Z_{k-1}},\ab{Z_k})\\
&&+C\theta\ab{z}^2\sqrt{\nu_1(k-3)\nu_1(k-2)\nu_1(k-1)\nu_1(k)}\\
&=& 1+\nu_1(k)z+C\theta\ab{z}^2\Big\{\nu_2(k)+\Expect
X_{k-1}X_k+\nu_1(k-1)\nu_1(k)+\sum_{l=0}^3\nu_1^2(k-l)\Big\}\\
&=&1+\nu_1(k)z+C\theta\ab{z}^2r_0(k),
\end{eqnarray*}
which proves (\ref{f1}).

\noindent The proof of  (\ref{f2}) is almost identical. We take longer expansion
in Lemma \ref{Hversion} and note that due to (\ref{bergident})
 \begin{equation*}
Z_k=X_kz+\theta X_k(X_k-1)\frac{\ab{z}^2}{2}.
\end{equation*}
Therefore,
\begin{eqnarray*}
\wE(Z_{k-1},Z_k)&=&\wE(X_{k-1}z+\theta\ab{z}^2X_{k-1}(X_{k-1}-1),Z_k)=z\wE
X_{k-1}Z_k\\
&&+C\theta\ab{z}^3 \wE^{+}(X_{k-1}(X_{k-1}-1),X_k)
=z^2\wE(X_{k-1},X_k)+C\theta\ab{z}^3\wE_2^{+}(X_{k-1},X_k).
\end{eqnarray*}
The other proofs are simple repetition of the given ones with the only
exception that results from previous steps are used. For example,
for the proof of (\ref{f3}), we apply Lemma \ref{Hversion} and get
\[
\ab{\vfi_k}=1+\Expect Z_k
 +\sum_{j=1}^{k-1}\frac{\ab{\wE(Z_j,\dots,Z_k)}}{\ab{\vfi_j\vfi_{j+1}\cdots\vfi_{k-1}}}\\
= 1+\Expect
Z_k+\sum_{j=k-2}^{k-1}+\sum_{j=k-6}^{k-3}+\sum_{j=1}^{k-7}.
\]
By (\ref{ez00}),
\[\sum_{j=1}^{k-7}\frac{\ab{\wE(Z_j,\dots,Z_k)}}{\ab{\vfi_j\vfi_{j+1}\cdots\vfi_{k-1}}}\leqslant
C\ab{z}^4\sqrt{\nu_1(k-1)\cdots\nu_1(k)}\leqslant
C\ab{z}^4\sum_{l=0}^7\nu_1^4(k-l)
\]
and by (\ref{shorting})
\[
\Ab{\sum_{j =k-6}^{k-3}}\leqslant
C\sum_{j=k-6}^{k-3}\wE^{+}(\ab{Z_j},\dots,\ab{Z_k})\leqslant
C\wE^{+}(\ab{Z_{k-3}},\dots,\ab{Z_k})\leqslant
C\ab{z}^4\wE^{+}(X_{k-3},\dots,X_k).
\]
\noindent For other summands, we apply Lemma \ref{factorial} and use the previous
estimates. \qed

Hereafter, the prime denotes the derivative with respect to $t$.
\begin{lemma}\label{prefis} Let condition (\ref{nu12}) hold. Then, for all
$t\in\RR$,
\begin{eqnarray*}
(\wE(Z_j,\dots,Z_k))'&=&\sum_{i=j}^k\wE(Z_j,\dots,Z_i',\dots,Z_k),\\
\ab{\wE(Z_j,\dots,Z_i',\dots,Z_k)}&\leqslant&
2^{[3(k-j)+1]/2}\ab{z}^{(k-j)/2}\prod_{l=j}^k\sqrt{\nu_1(l)}.
\end{eqnarray*}
\end{lemma}
\noindent \Proof The first identity was proved in \citep{H82}. Applying
(\ref{ezkj}) we obtain
\[
\ab{\wE(Z_j,\dots,Z_i',\dots,Z_k)}\leqslant
2^{k-j}\sqrt{\Expect\ab{Z_i'}^2}\prod_{l\ne
i}^k\sqrt{\Expect\ab{Z_l}^2}.
\]
Due to assumption (\ref{nu12}), $\nu_2(l)\leqslant \nu_1(l)$. Therefore,

\[
\Expect\ab{Z_l'}^2=\Expect\ab{\ii\ee^{\ii t X_l}X_l}^2=\Expect
X_l^2=\Expect X_l(X_l-1+1)=\nu_2(l)+\nu_1(l)\leqslant 2\nu_1(l).
\]
Combining the last estimate with $\Expect\ab{Z_l}^2\leqslant
2\Expect\ab{Z_l}\leqslant 2\ab{z}\nu_1(l)$, the proof follows. \qed

\begin{lemma}\label{fis} Let condition (\ref{nu12}) be satisfied, $k=1,\dots,n$ and $\vfi_k$ be defined as in Lemma
  \ref{Hversion}. Then, for all $t\in\RR$,
\begin{eqnarray}
\vfi_k'&=&33\theta[\nu_1(k)+\nu_1(k-1)],\label{fis0}\\
\vfi_k'&=&\nu_1(k)z'+\theta C\ab{z}(r_0(k)+\wE^{+}(X_{k-2},X_{k-1})),\label{fis1}
\\
\vfi_k'&=&\nu_1(k)z'+\gamma_2(k)(z^2)'+\theta
C\ab{z}^2\Big(r_1(k)+[\nu_1(k-2)+\nu_1(k)]\Expect X_{k-1}X_k\nonumber\\
&&+\wE^{+}(X_{k-4},X_{k-3},X_{k-2})+\wE^{+}(X_{k-3},X_{k-2},X_{k-1})\Big),\label{fis2}\\
\vfi_k'&=&\nu_1(k)z'+\gamma_2(k)(z^2)'+\gamma_3(k)(z^3)'+\theta
C\ab{z}^3\Big(r_2(k)+\wE^{+}(X_{k-4},\dots,X_{k-1})\nonumber\\
&&+\wE^{+}(X_{k-5},\dots,X_{k-2})\Big).\label{fis3}
\end{eqnarray}
\end{lemma}

\noindent \Proof Note that
\[
\bigg(\frac{\wE(Z_j,\dots,Z_k)}{\vfi_j\cdots\vfi_{k-1}} \bigg)'=
\frac{(\wE(Z_j,\dots,Z_k))'}{\vfi_j\dots\vfi_k}
-\frac{\wE(Z_j,\dots,Z_k)}{\vfi_j\dots\vfi_k}\sum_{m=j}^{k-1}\frac{\vfi_m'}{\vfi_m}.
\]
Now the proof is just a repetition of the proof of Lemma \ref{f}.
For example, (\ref{fis0}) is easily verifiable for $k=0,1$. Let us
assume that it holds for $1,2,\dots,k-1$.
 From Lemmas \ref{Hversion} and \ref{factorial} and equation
(\ref{sh1}),  we get
\begin{eqnarray*}
\ab{\vfi_k'}&\leqslant&
\nu_1(k)+\sum_{j=1}^{k-1}\frac{\ab{(\wE(Z_j,\dots,Z_k))'}}{\ab{\vfi_j\cdots\vfi_{k-1}}}+
\sum_{j=1}^{k-1}\frac{\ab{\wE(Z_j,\dots,Z_k)}}{\ab{\vfi_j\cdots\vfi_{k-1}}}\sum_{m=j}^{k-1}\frac{\ab{\vfi_m'}}{\ab{\vfi_m}}
\\
&\leqslant&
\nu_1(k)+\sum_{j=1}^{k-1}\bigg(\frac{10}{9}\bigg)^{k-j}
\sum_{i=j}^k\ab{\wE(Z_j,\dots,Z_i',\dots,Z_k)}\\
&&+ \sum_{j=1}^{k-1}\bigg(\frac{10}{9}\bigg)^{k-j}
\ab{\wE(Z_j,\dots,Z_k)}(k-j)33\cdot0.02\cdot\bigg(\frac{10}{9}\bigg).
\end{eqnarray*}
By Lemma \ref{prefis},
\[
\ab{\wE(Z_j,\dots,Z_i',\dots,Z_k)}\leqslant
[\nu_1(k-1)+\nu_1(k)](0.04)^{k-j}\frac{10}{\sqrt{2}}.
\]
Combining the last two estimates and (\ref{ez00}), the
proof of (\ref{fis0}) is completed.

\noindent We omit the proofs of remaining expansions and note
only that
\[
(\ee^{\ii t X}-1)'=\ii X\ee^{\ii t X}=\ii \eit X\ee^{\ii t(X-1)}=
z'X\Big(1+(X-1)z+\theta\frac{(X-1)(X-2)}{2}\ab{z}^2\Big),
\]
 due to Bergstr\"om's identity. \hfill{\qed}

Let, for $j=1,\dots,n$ and $l=2,\dots,n$,
\begin{equation*}
 g_j(t)=\Exponent{\nu_1(j)(\ee^{\ii
t}-1)+\bigg(\frac{\nu_2(j)-\nu_1^2(j)}{2}+\wE(X_{j-1},X_j)
\bigg)(\ee^{\ii t}-1)^2}.\label{gj}
 \end{equation*}

\begin{lemma}\label{gkisv} Let conditions in (\ref{nu12}) be
satisfied, $k=1,2,\dots,n$. Then, for all $t\in\RR$,
\begin{eqnarray}
g_k&=&1+C\theta\ab{z}[\nu_1(k-1)+\nu_1(k)],\label{g0}\\
g_k'&=&C\theta[\nu_1(k-1)+\nu_1(k)],\label{gis0}
\\
g_k&=&1+\nu_1(k)z+\gamma_2(k)z^2+C\theta\ab{z}^3\big\{\nu_1^3(k-1)+\nu_1^3(k)+\nu_1(k)\nu_2(k)\nonumber\\
&&+[\nu_1(k-1)+\nu_1(k)]\Expect X_{k-1}X_k\big\},\label{g1}\\
g_k'&=&\nu_1(k)z'+\gamma_2(k)(z^2)'+C\theta\ab{z}^2\big\{\nu_1^3(k-1)+\nu_1^3(k)+\nu_1(k)\nu_2(k)\nonumber\\
&&+[\nu_1(k-1)+\nu_1(k)]\Expect X_{k-1}X_k\big\},\label{gis1}\\
g_k&=&1+\nu_1(k)z+\gamma_2(k)z^2+\tilde\gamma_3(k)z^3\nonumber\\
&&+C\theta\ab{z}^4\big\{\nu_1^4(k-1)+\nu_1^4(k)+\nu_2^2(k)+(\Expect X_{k-1}X_k)^2\big\},\label{g2}\\
g_k'&=&\nu_1(k)z'+\gamma_2(k)(z^2)'+\tilde\gamma_3(k)(z^3)'\nonumber\\
&&+C\theta\ab{z}^3\big\{\nu_1^4(k-1)+\nu_1^4(k)+\nu_2^2(k)+(\Expect X_{k-1}X_k)^2\big\},\label{gis2}\\
\ab{g_k}&\leqslant&\exponent{-\lambda_k\sin^2(t/2)}.\label{glambda}
\end{eqnarray}
Here $\lambda_k$ is as in Lemma \ref{lambdak} and
\[\tilde \gamma_3(k)=
\frac{\nu_1(k)\nu_2(k)-\nu_1^3(k)}{2}+\nu_1(k)\wE(X_{k-1},X_k)+\frac{\nu_1^3(k)}{6}.
\]
\end{lemma}

\noindent  \Proof For any complex number $b$, we have
\[
\ee^b=1+b+\frac{b^2}{2}+\dots+\frac{b^s}{s!}+\theta\frac{\ab{b}^{s+1}}{(s+1)!}\ee^{\ab{b}}.
\]
 Due to (\ref{nu12}), $\nu_2(j)\leqslant \nu_1(j)$. Therefore,
\begin{eqnarray}
\Expect X_{j-1}X_j&\leqslant& \sqrt{\Expect X_{j-1}^2\Expect
X_j^2}\leqslant
\sqrt{[\nu_2(j-1)+\nu_1(j-1)][\nu_2(j)+\nu_1(j)]}\nonumber\\
&\leqslant& 2\sqrt{\nu_1(j-1)\nu(j)}\leqslant
2[\nu_1(j-1)+\nu_1(j)]. \label{ex1x2}
\end{eqnarray}
Therefore, the exponent of $g_k$ is bounded by some absolute
constant $C$ and (\ref{g0}) and (\ref{gis0}) easily follow. We
have
\begin{eqnarray*}
g_k&=&1+\nu_1(k)z+\gamma_2(k)z^2+C\theta\big\{\nu_1^3(k)+\nu_2^2(k)+\nu_1(k)\nu_2(k)\\
&&+\nu_1(k)\wE^{+}(X_{k-1},X_k)+\nu_1(k-1)\nu_1^2(k)+(\wE^{+}(X_{k-1},X_k))^2
\big\}.
\end{eqnarray*}
Moreover,
\[\nu_2^2(k)\leqslant \nu_1(k)\nu_2(k),\qquad
\nu_1(k-1)\nu_1^2(k)\leqslant \nu_1^3(k-1)+\nu_1^3(k)\] and
\[
(\wE^{+}(X_{k-1},X_k)\leqslant 2(\Expect
X_{k-1}X_k)^2+2\nu_1^2(k-1)\nu_1^2(k)\leqslant
2[\nu_1(k-1)+\nu_1(k)]\Expect
X_{k-1}X_k+2\nu_1^3(k-1)+2\nu_1^3(k).
\]
Thus, (\ref{g1}) easily follows. The estimates (\ref{gis1}) --
(\ref{gis2}) are proved similarly.

\noindent  For the proof of
(\ref{glambda}), note that
 \[\wE^{+}(X_{k-1},X_k)\leqslant \Expect
 X_{k-1}X_k+0.01\nu_1(k),\quad \nu_1^2(k)\leqslant 0.01\nu_1(k)\]
 and
 \begin{eqnarray*}
 \ab{g_k}&\leqslant&\Exponent{-2\nu_1(k)\sin^2(t/2)+2[\nu_2(k)+\nu_1^2(k)+2\wE^{+}(X_{k-1},X_k)]\sin^2(t/2)}\\
 &\leqslant&\Exponent{-1.92\nu_1(k)\sin^2(t/2)+2\nu_2(k)\sin^2(t/2)+4\Expect
 X_{k-1}X_k\sin^2(t/2)},
 \end{eqnarray*}
which completes the proof. \qed

For asymptotic expansions, we need a few smoothing estimates.

\begin{lemma} \label{smoothing} Let conditions (\ref{nu12}) and (\ref{lambda}) be
satisfied, $0\leqslant \alpha\leqslant 1$, and $M$ be any
finite (signed) measure. Then
\begin{equation*}
\norm{M\exponent{\Gamma_1U+\alpha\Gamma_2U^2}}\leqslant
C\norm{M\exponent{0.9\lambda U}}.
\end{equation*}
\end{lemma}

\noindent \Proof Due to (\ref{nu12}) and (\ref{lambda}), we have
\begin{equation*}
\Gamma_1-3.1\ab{\Gamma_2}\geqslant\Gamma_1-1.55\sum_{k=1}^n\nu_2(k)-0.0155\Gamma_1-3.1\sum_{k=1}^n\Expect
X_{k-1}X_k-0.031\Gamma_1\geqslant 0.9\lambda. \label{sm1}
\end{equation*}
Thus, \begin{eqnarray*}
\norm{M\exponent{\Gamma_1U+\alpha\Gamma_2U^2}}&\leqslant&
\norm{M\exponent{(\Gamma_1-3.1\ab{\Gamma_2})U}}
\norm{\exponent{3.1\ab{\Gamma_2}U+\alpha\Gamma_2U^2}}\\
&\leqslant&\norm{M\exponent{0.9\lambda U}}
\norm{\exponent{3.1\ab{\Gamma_2}U+\alpha\Gamma_2U^2}}.
\end{eqnarray*}
It remains to prove that
the second exponent measure is bounded by some absolute constant.
Note that the total variation of \emph{any distribution} equals unity.
Therefore, by Lemma \ref{a10}
\begin{eqnarray*}
\lefteqn{\Norm{\exponent{3.1\ab{\Gamma_2}U+\alpha\Gamma_2U^2}}=
 \Norm{\exponent{3.1\ab{\Gamma_2}U}\Big(\dirac+\sum_{m=1}^\infty\frac{(\alpha\Gamma_2U^2)^m}{m!}\Big)}}\hskip 0.5cm\\
 &\leqslant&1+\sum_{m=1}^\infty\frac{\ab{\Gamma_2}^m}{m!}\norm{U^2\exponent{3.1\ab{\Gamma_2}U/m}}^m\leqslant
 1+\sum_{m=1}^\infty\frac{\ab{\Gamma_2}^m}{m^m\ee^{-m}\sqrt{2\pi m}}\bigg(\frac{3m}{3.1\ab{\Gamma_2}\ee}\bigg)^m\leqslant
 C.
\end{eqnarray*}
Combining both inequalities given above, we complete the proof of
the lemma. \qed

\begin{lemma} \label{NBsmoothing} Let conditions (\ref{nu12}) and (\ref{3ab})
be satisfied. Then
\begin{eqnarray}
\NB(r,\qubar)&=&\Exponent{\Gamma_1U+\Gamma_2U^2+\frac{4\Gamma_2^2}{3\Gamma_1}U^3+
\frac{2\Gamma_2^3}{\Gamma_1^2}U^4\Theta\frac{1}{0.7}}
=\Exponent{\Gamma_1U+\Gamma_2U^2+\frac{4\Gamma_2^2}{3\Gamma_1}U^3\Theta\frac{1}{0.7}}\nonumber\\
&=&\Exponent{\Gamma_1U+\Gamma_2U^2\Theta\frac{1}{0.7}}=
\Exponent{\Gamma_1U+\Gamma_1U^2\Theta\frac{3}{28}}=\Exponent{0.5\Gamma_1U}\Theta
C.\label{NBsmu}
\end{eqnarray}
\end{lemma}

\noindent \Proof  Due to (\ref{3ab}),
\[
\Gamma_2=\frac{1}{2}\sum_{k=1}^n(\nu_2^2(k)-\nu_1^2(k))+\sum_{k=1}^n
Cov(X_{k-1},X_k)\leqslant
\frac{1}{2}\sum_{k=1}^n\nu_2(k)+\sum_{k=1}^n\ab{Cov(X_{k-1},X_k)}\leqslant
\frac{3}{40}\Gamma_1.
\]
Therefore,
\[
\frac{1-\qubar}{\qubar}=\frac{2\Gamma_2}{\Gamma_1}\leqslant 0.15,
\quad \bigg(\frac{1-\qubar}{\qubar}\bigg)\norm{U}\leqslant
0.15(\norm{\dirac_1}+\norm{\dirac})\leqslant 0.3.
\]
Consequently, from (\ref{nbinom}),
\begin{eqnarray*}
\NB(r,\qubar)&=&\Exponent{\sum_{j=1}^\infty
\frac{r}{j}\bigg(\frac{1-\qubar}{\qubar} \bigg)^j U^j}\\
&=&
\Exponent{\Gamma_1U+r\bigg(\frac{1-\qubar}{\qubar}\bigg)^2\frac{U^2}{2}+
r\bigg(\frac{1-\qubar}{\qubar}\bigg)^3\frac{U^3}{3}+r\bigg(\frac{1-\qubar}{\qubar}\bigg)^4\frac{U^4}{4}\Theta\frac{1}{0.7}}\\
&=&
\Exponent{\Gamma_1U+r\bigg(\frac{1-\qubar}{\qubar}\bigg)^2\frac{U^2}{2}+
r\bigg(\frac{1-\qubar}{\qubar}\bigg)^3\frac{U^3}{3}\Theta\frac{1}{0.7}}\\
&=&
\Exponent{\Gamma_1U+r\bigg(\frac{1-\qubar}{\qubar}\bigg)^2\frac{U^2}{2}\Theta\frac{1}{0.7}}.
\end{eqnarray*}
Recalling that $r(1-\qubar)/\qubar=\Gamma_1$, we obtain all
equalities except the last one. The last equality is equivalent to
\[
\Norm{\Exponent{0.5\Gamma_1U+\Gamma_1U^2\Theta\frac{3}{28}}}\leqslant
C
\]
which is proved similarly to Lemma \ref{smoothing}. \qed

\begin{lemma} \label{BINsmooth}Let  conditions (\ref{nu12}) and (\ref{3ab}) be satisfied. Then
\begin{eqnarray*}
\BI(N,\pbar)&=&\Exponent{-N\sum_{j=1}^\infty
\frac{(-\pbar U)^j}{j}}\nonumber\\
&=& \Exponent{\Gamma_1U
+\Gamma_2U^2+U^2\theta\frac{50\Gamma_2^2\epsilon}{21\Gamma_1^2}
+\frac{N\pbar^3
U^3}{3}+\frac{N\pbar^4U^4}{4}\Theta\frac{5}{3}}\nonumber\\
&=&
\Exponent{\Gamma_1U+\Gamma_2U^2+U^2\theta\frac{50\Gamma_2^2\epsilon}{21\Gamma_1^2}+\frac{N\pbar^3
U^3}{3}\Theta\frac{5}{3}}\nonumber\\
&=& \Exponent{\Gamma_1U+\frac{N\pbar^2 U^2}{2}\Theta\frac{5}{3}}
=\Exponent{\Gamma_1U+\Gamma_1 U^2\Theta\frac{1}{6}}
=\Exponent{0.5\Gamma_1U}\Theta C.\label{BIsmu}
\end{eqnarray*}
\end{lemma}

\noindent \Proof  Due to (\ref{3ab}),
\begin{equation*}
\ab{\Gamma_2}\leqslant\frac{1}{2}\sum_{k=1}^n(\nu_2(k)+0.01\nu_1(k))+\sum_{k=1}^n
\ab{Cov(X_{k-1},X_k)}\leqslant
\Gamma_1(0.025+0.005+0.05)=0.08\Gamma_1.
\label{bgam2}\end{equation*}
  Therefore,
\begin{equation}
\pbar=\frac{\Gamma_1}{\tilde
N-\epsilon}\leqslant\frac{\Gamma_1}{\tilde
N-1}\leqslant\frac{2\ab{\Gamma_2}}{\Gamma_1-2\ab{\Gamma_2}}\leqslant\frac{50\ab{\Gamma_2}}{21\Gamma_1}<\frac{1}{5}.
\label{barp}\end{equation}
 and
\[
\frac{\epsilon}{\tilde N}\leqslant\frac{1}{\tilde
N}=\frac{2\ab{\Gamma_2}}{\Gamma_1^2}\leqslant\frac{2\ab{\Gamma_2}}{\Gamma_1}\leqslant
0.16.
\]
Consequently,
\[
N\pbar^2=2\ab{\Gamma_2}\frac{\tilde N}{N}=
2\ab{\Gamma_2}\frac{1}{1-\epsilon/\tilde N}=
2\ab{\Gamma_2}\bigg(1+\frac{\epsilon}{\tilde
N}\theta\frac{100}{84}\bigg)
\]
and
\begin{equation}
-\frac{N\pbar^2}{2}=\Gamma_2+\theta\frac{50\Gamma_2^2}{21\Gamma_1^2}\epsilon.
\label{bin2}
\end{equation}
Taking into account (\ref{barp}), we prove
\begin{eqnarray*}
\BI(N,\pbar)&=&\Exponent{-N\sum_{j=1}^\infty
\frac{(-\pbar U)^j}{j}}\\
&=& \Exponent{\Gamma_1U-\frac{N(\pbar U)^2}{2}+\frac{N(\pbar
U)^3}{3}+\frac{N\pbar^4U^4}{4}\Theta\frac{5}{3}}\\
&=& \Exponent{\Gamma_1U-\frac{N(\pbar U)^2}{2}+\frac{N(\pbar
U)^3}{3}\Theta\frac{5}{3}}\\
&=& \Exponent{\Gamma_1U+\frac{N(\pbar U)^2}{2}\Theta\frac{5}{3}}
=\Exponent{\Gamma_1U+\Gamma_1 U^2\Theta\frac{1}{6}}.
\end{eqnarray*}
Combining (\ref{bin2}) with the last expansions, we  obtain all
equalities except the last one whose proof is similar to that of Lemma
\ref{smoothing}. \qed

\section{Proofs }

\noindent \textbf{Proof of Theorem \ref{G2T}}. Let $\w M(t)=J_1+J_2$, where
\begin{eqnarray*}
J_1&=&\prod_{j=1}^n\vfi_j-\prod_{j=1}^ng_j-\sum_{m=1}^{n}(\vfi_m-g_m)\prod_{j\ne m}^n g_m,\\
J_2&=&\prod_{j=1}^ng_j+\sum_{m=1}^n(\vfi_m-g_m)\prod_{j\ne
m}g_j-\prod_{j=1}^ng_j(1+\Gamma_3z^3).
\end{eqnarray*}

We estimate $J_1$ and $J_2$ separately. Further we frequently apply the following  estimate
 \begin{eqnarray}
\prod_{j=1,j\ne m,l}^n\exponent{-\lambda_j\sin^2(t/2)}&\leqslant&
\exponent{-1.3\lambda\sin^2(t/2)}\exponent{(\lambda_m+\lambda_l)\sin^2(t/2)}\nonumber\\
&\leqslant& C\exponent{-1.3\lambda\sin^2(t/2)},\label{naujas2}
\end{eqnarray}
which is valid for any $m,l\in\{1,2,\dots,n\}$, since
all $\lambda_j\leqslant C$.

 Applying the generalized Bergstr\"om identity from \citep{Ce98},  (\ref{naujas2}), Lemmas
 \ref{lambdak},  \ref{f} , \ref{gkisv} and \ref{simple}, we obtain
\begin{eqnarray*}
\ab{J_1}&=&\bigg\vert\sum_{l=2}^n(\vfi_l-g_l)\prod_{j=l+1}^n\vfi_j\sum_{m=1}^{l-1}(\vfi_m-g_m)\prod_{j=1,j\ne m}^{l-1}g_j\bigg\vert\nonumber\\
&\leqslant&C\sum_{l=2}^n\ab{\vfi_l-g_l}\sum_{m=1}^{l-1}\ab{\vfi_m-g_m}\prod_{j=1,j\ne m,l}^n\exponent{-1.3\lambda_j\sin^2(t/2)}\nonumber\\
&\leqslant& C\exponent{-1.3\lambda
\sin^2(t/2)}\bigg(\sum_{k=1}^n\ab{\vfi_k-g_k}\bigg)^{2}
\leqslant  C\exponent{-1.3\lambda
\sin^2(t/2)}R_1^2\ab{z}^6\nonumber\\
&\leqslant& C\exponent{-\lambda
\sin^2(t/2)}R_1^2\min(1,\lambda^{-3}). \label{naujas3}
\end{eqnarray*}

Similarly, taking into account
 (\ref{f3}),(\ref{g0}), (\ref{g2}) and
(\ref{glambda}), we get
\begin{eqnarray*}
\ab{J_2}
&\leqslant&\Ab{\prod_{j=1}^ng_j(1+\Gamma_3
z^3)-\prod_{j=1}^n g_j-\prod_{j=1}^ng_j\sum_{m=1}^n(\vfi_m-g_m)}\\
&&+
\Ab{\sum_{m=1}^n(\vfi_m-g_m)\bigg(\prod_{j=1}^n g_j-\prod_{j\ne
m}^ng_j\bigg)}\\
&=&\Ab{\prod_{j=1}^ng_j\Big(\sum_{m=1}^n(\vfi_m-g_m)-\Gamma_3z^3\Big)}+
\Ab{\sum_{m=1}^n(\vfi_m-g_m)\prod_{j\ne m} g_j(g_m-1)}\\
&\leqslant& C
R_2\ab{z}^4\exponent{-1.3\lambda\sin^2(t/2)}\leqslant C\exponent{-\lambda\sin^2(t/2)}R_2\min(1,\lambda^{-2}).
\end{eqnarray*}
Therefore,
\begin{equation}
\ab{\w M(t)}\leqslant C\exponent{-\lambda
\sin^2(t/2)}(R_1^2\min(1,\lambda^{-3})+R_2\min(1,\lambda^{-2})).
\label{naujas4}
\end{equation}

Let $\tilde\vfi_k=\vfi_k\exponent{-\ii\nu_1(k)t}$, $\tilde g_k=
g_k\exponent{-\ii\nu_1(k)t}$ . Observe that
$
\ab{\tilde\vfi_l'-\tilde g_l'}\leqslant
C(\ab{\vfi_l'-g_l'}+\nu_1(k)\ab{\vfi_l-g_l})$.
Moreover, taking into account (\ref{fis1}), (\ref{fis2}) and (\ref{ex1x2}),
we get
\[
\ab{\tilde\vfi_l'}\leqslant\ab{\vfi_l'-\nu_1(l)z'}+\nu_1(l)\ab{\eit-\vfi_l}\leqslant
\ab{\vfi_l'-\nu_1(l)z'}+\nu_1(l)\ab{z}+\nu_1(l)\ab{1-\vfi_l}\leqslant
C\ab{z}\sum_{j=0}^3\nu_1(l-j)
\]
and similar estimate holds for $\ab{\tilde g_l'}$.

Taking into account (\ref{naujas2}), Lemmas
 \ref{lambdak},  \ref{f}, \ref{fis} and \ref{gkisv} we prove that

 \begin{eqnarray}
\Ab{(\ee^{-\ii\Gamma_1 t}J_1)'}&\leqslant&\sum_{l=2}^n\ab{\tilde\vfi_l'-\tilde g_l'}\prod_{j=l+1}^n\ab{\tilde\vfi_l}\sum_{m=1}^{l-1}
\ab{\tilde\vfi_m-\tilde g_m}\prod_{j=1,j\ne m}^{l-1}\ab{\tilde g_j}\nonumber\\
&&+\sum_{l=2}^n\ab{\tilde\vfi_l-\tilde g_l}\sum_{j=l+1}^n\ab{\tilde\vfi_j'}\prod_{i=l+1,i\ne j}^n\ab{\tilde\vfi_j}\sum_{m=1}^{l-1}
\ab{\tilde\vfi_m-\tilde g_m}\prod_{j=1,j\ne m}^{l-1}\ab{\tilde g_j}\nonumber\\
&&+\sum_{l=2}^n\ab{\tilde\vfi_l-\tilde g_l}\prod_{j=l+1}^n\ab{\tilde\vfi_j}\sum_{m=1}^{l-1}\ab{\tilde\vfi_m'-\tilde g_m'}\prod_{j=1,j\ne m}^{l-1}\ab{\tilde g_j}\nonumber\\
&&+\sum_{l=2}^n\ab{\tilde\vfi_l-\tilde g_l}\prod_{j=l+1}^n\ab{\tilde\vfi_j}\sum_{m=1}^{l-1}\ab{\tilde\vfi_m-\tilde g_m}
\sum_{j=1,j\ne m}^{l-1}\ab{\tilde g_j'}\prod_{k=1,k\ne m,j}^{l-1}\ab{\tilde g_k}\nonumber\\
&\leqslant& C\exponent{-1.3\lambda\sin^2(t/2)}\bigg(\sum_{l=2}^n\ab{\tilde\vfi_l'-\tilde g_l'}\sum_{m=1}^{l-1}\ab{\tilde\vfi_m-\tilde g_m}\nonumber\\
&&+\bigg(\sum_{l=1}^n\ab{\tilde\vfi_l-\tilde g_l}\bigg)^2\sum_{j=1}^n(\ab{\tilde\vfi_j'}+\ab{\tilde g_j'})\bigg)\nonumber\\
&\leqslant&C\exponent{-1.3\lambda\sin^2(t/2)}(R_1^2\ab{z}^5(1+\Gamma_1\ab{z}^2)+R_1^2\ab{z}^7\Gamma_1)\nonumber\\
&\leqslant&C\exponent{-\lambda\sin^2(t/2)}(1+\Gamma_1\min(1,\lambda^{-1}))R_1^2\min(1,\lambda^{-5/2}).\label{naujas5}
 \end{eqnarray}
Similarly
 \begin{eqnarray*}
\ab{(\ee^{-\ii\Gamma_1 t}J_2)'}
&\leqslant& \Ab{\Big(\ee^{-\ii
t\Gamma_1}\prod_{j=1}^ng_j\sum_{m=1}^n(\vfi_m-g_m-\gamma_3(m)z^3)\Big)'}\\
&&+
\Ab{\Big(\ee^{-\ii t\Gamma_1}\sum_{m=1}^n(\vfi_m-g_m)\prod_{j\ne
m}g_j(g_m-1)\Big)'}\\
&\leqslant&\Ab{\Big(\prod_{j=1}^n\tilde
g_j\Big)'\sum_{m=1}^n(\vfi_m-g_m-\gamma_3(m)z^3)}+\Ab{\prod_{j=1}^n\tilde
g_j\sum_{m=1}^n(\vfi_m'-g_m'-\gamma_3(m)(z^3)')}\\
&&+\Ab{\sum_{m=1}^n(\vfi_m'-g_m')\prod_{j\ne m}\tilde g_j(\tilde
g_m-\ee^{-\ii
t\nu_1(m)})}+\sum_{m=1}^n\ab{\vfi_m-g_m}\ab{g_m-1}\Ab{\Big(\prod_{j\ne
m}\tilde g_j\Big)'}\\
&&+\sum_{m=1}^n\ab{\vfi_m-g_m}\Ab{\prod_{k\ne m}\tilde g_j}[\ab{\tilde
g_m'}+\nu_1(m)].
 \end{eqnarray*}

Applying Lemmas \ref{fis}, \ref{gkisv} and \ref{simple}, it is not
difficult to prove that the derivative given above is less than
$C\ab{z}^5\Gamma_1R_2\exponent{-1.3\lambda\sin^2(t/2)}$. Combining
this estimate with (\ref{naujas5}) we obtain
\[\ab{(\ee^{-\ii\Gamma_1 t}\w M(t))'}\leqslant C\exponent{-\lambda\sin^2(t/2)}(1+\Gamma_1\min(1,\lambda^{-1})(R_1^2\min(1,\lambda^{-5/2})+R_2\min(1,\lambda^{-3/2})).\]

 For the proof of (\ref{GVA}), we
use (\ref{naujas4}),  (\ref{TVAP}) with
$v=\Gamma_1$ and $u=\max(1,\Gamma_1)$. For the proof of (\ref{GV}) we use identity
\begin{equation}
\prod_{j=1}^n\vfi_j-\prod_{j=1}^n g_j=\sum_{j=1}^n(\vfi_j-g_j)\prod_{l=j+1}^n\vfi_l\prod_{l=1}^{j-1}g_l.\label{naujas6}
\end{equation}
The rest of the proof is  very similar to the proof of (\ref{GVA}) and, therefore, omitted.
 \qed

\noindent \textbf{Proof of Theorem \ref{PoissonT}}. For the proof of (\ref{PV}) we use (\ref{naujas6}) with $g_j$ replaced by $\exponent{\nu_1(j)z}$. Now the proof is very similar to the proofs of (\ref{GVA}) and (\ref{GV}) and, therefore,  omitted.
Applying Lemma \ref{smoothing} and using the following identity
\begin{equation}
\ee^b-1-b=b^2\int_0^1(1-\tau)\ee^{\tau b}\dd\tau \label{tau},
\end{equation}
 we get
\begin{eqnarray*}
\lefteqn{ \Norm{\G-\Pois(\Gamma_1)(\dirac+\Gamma_2U^2)}=
\Norm{\exponent{\Gamma_1
U}\int_0^1(1-\tau)(\gamma_2U^2)^2\exponent{\tau\Gamma_2U^2}\dd
\tau}}\hskip 0.5cm\\
&\leqslant&\int_0^1\norm{\Gamma_2^2U^4\exponent{\Gamma_1U+\tau\Gamma_2U^2}}
\dd\tau\leqslant C\ab{\Gamma_2}^2\norm{U^4\exponent{0.9\lambda
U}}\leqslant CR_0^2\min(1,\lambda^{-2}).
\end{eqnarray*}
Combining this estimate with Bergstr\"om expansion ($s=1$) for $\G$,
we prove (\ref{PVA}).   \qed

\noindent \textbf{Proof of Theorem \ref{NBtheorem}.} Applying (\ref{NBsmu})
and Lemma \ref{a10}, we obtain
\begin{eqnarray*}
\Norm{\G-\NB(r,\qubar)}&=&\Norm{\G-\G\Exponent{\frac{4\Gamma_2^2}{\Gamma_1}U^3\Theta\frac{1}{0.7}}}=
\Norm{\G\int_0^1\Big(\Exponent{\tau\frac{4\Gamma_2^2}{\Gamma_1}U^3\Theta\frac{1}{0.7}}\Big)'\dd\tau}\\
&\leqslant&C\int_0^1\frac{\Gamma_2^2}{\Gamma_1}\Norm{U^3\Exponent{\Gamma_1U+\Gamma_2U^2+\tau\frac{4\Gamma_2^2}{\Gamma_1}U^3\Theta\frac{1}{0.7}}}
\dd\tau\\
&\leqslant&
C\frac{\Gamma_2^2}{\Gamma_1}\Norm{U^3\exponent{0.5\Gamma_1U}}\leqslant
C\frac{\Gamma_2^2}{\Gamma_1}\min(1,\Gamma_1^{-3/2}).
\end{eqnarray*}
Combining the last estimate with (\ref{GV}), we prove (\ref{NBV}).

Let
\[
M_1:=\frac{4\Gamma_2^2}{3\Gamma_1}U^3,\quad M_2:=
\frac{2\Gamma_2^3}{\Gamma_1^2}U^4\Theta\frac{1}{0.7},\quad M_3:=
\Gamma_3U^3-M_1.
\]
 Then by Lemmas
\ref{NBsmoothing} and \ref{a10} and using equation (\ref{tau}),
\begin{eqnarray*}
\NB(r,\qubar)&=&\G\exponent{M_1+M_2}\\
&=&\G\Big(\dirac+M_1+M_1^2\int_0^1(1-\tau)\exponent{\tau
M_1}\dd\tau \Big)\Big(\dirac+M_2\int_0^1\exponent{xM_2}\dd
x \Big)\\
&=&\G(\dirac+M_1)+M_1^2\int_0^1(1-\tau)\G\exponent{\tau
M_1}\dd\tau\\
&&+ \int_0^1\int_0^1M_2(\dirac+M_1+M_1^2(1-\tau))\G\exponent{\tau
M_1+xM_2}\dd\tau\dd x\\
&=&\G(\dirac+M_1)+\exponent{0.5\Gamma_1U}(M_1^2\Theta
C+[M_2+M_1M_2]\Theta C +M_1^2M_2\Theta C)\\
&=&\G(\dirac+M_1)+\exponent{0.25\Gamma_1U}\Gamma_2^3\Gamma_1^{-2}U^4\Theta
C.
\end{eqnarray*}
By the triangle inequality,
\begin{eqnarray*}
\lefteqn{\norm{F_n-\NB(r,\qubar)(\dirac+M_3)}}\hskip 1cm\nonumber\\
&\leqslant&\norm{F_n-\G(\dirac+\Gamma_3U^3)}+\norm{\G(\dirac+\Gamma_3U^3)-
\G(\dirac+M_1)(\dirac+M_3)}\nonumber\\
&&+C\norm{\exponent{0.25\Gamma_1U}\Gamma_2^3\Gamma_1^{-1}U^4(\dirac+M_3)}=:J_{31}+J_{32}+J_{33}.
\label{sa1}
\end{eqnarray*}
By Lemmas \ref{smoothing} and  \ref{a10},
\[
J_{32}\leqslant C\norm{\exponent{0.9\lambda
U}\Gamma_2^2\Gamma_1^{-1}(\Gamma_3-4\Gamma_2^2(3\Gamma_1)^{-1})U^6}\leqslant
\Gamma_2^2\Gamma_1^{-1}\ab{\Gamma_3-4\Gamma_2^2(3\Gamma_1)^{-1}}\min(1,\Gamma_1^{-3}).
\]
Similarly
\begin{eqnarray*}
J_{33}&\leqslant&
C\norm{\exponent{0.25\Gamma_1U}\Gamma_2^3\Gamma_1^{-2}U^4}+C\norm{
\exponent{0.25\Gamma_1U}\Gamma_2^3\Gamma_1^{-2}
(\Gamma_3-4\Gamma_2^2(3\Gamma_1)^{-1})U^7}\\
&\leqslant&
C\Gamma_2^3\Gamma_1^{-2}\min(1,\Gamma_1^{-2})+C\Gamma_2^3\Gamma_1^{-2}
\ab{\Gamma_3-4\Gamma_2^2(3\Gamma_1)^{-1}}\min(1,\Gamma_1^{-7/2}).
\end{eqnarray*}
Combining the last two estimates and applying (\ref{GVA}) for
 $J_{31}$, we prove (\ref{NBVA}).  \qed

\noindent \textbf{Proof of Theorem \ref{Bintheorem}.} Let
\begin{equation} \label{neqn4}
\tilde M_1:=\frac{N\pbar^3U^3}{3},\quad \tilde M_2:=
\frac{N\pbar^4U^4}{4}\Theta\frac{5}{3}+U^2\theta\frac{50\Gamma_2^2\epsilon}{21\Gamma_1^2},\quad
\tilde M_3:= \Gamma_3U^3-\tilde M_1.
\end{equation}
Since the proof is almost identical to that  of Theorem
\ref{NBtheorem}, it is omitted.  \qed

\noindent \textbf{Proof of Theorem \ref{NBASHARP}.} Let $\tilde
M_3$ be defined as (\ref{neqn4}). Observe that
\begin{eqnarray*}
&&\nu_1(k)=p^2, \quad \nu_2(k)=\nu_3(k)=0,\quad \Expect
X_{k-1}X_k\leqslant Cp^3,\quad\Expect X_{k-2}X_{k-1}X_k\leqslant
Cp^4,\\
&&\Expect X_{k-3}\cdots X_k\leqslant Cp^5,\quad \Gamma_2\leqslant
Cnp^3,\quad\Gamma_3\leqslant Cnp^4,\quad R_1\leqslant Cnp^4,\quad
R_2\leqslant Cnp^5.
\end{eqnarray*}
and
\[\tilde M_3=-\frac{np^4}{3}U^3+U^3\theta Cnp^5.\]
From Lemmas \ref{NBsmoothing} and \ref{a10}, we have
\begin{eqnarray}
\Norm{\big(\NB(r,\qubar)-\exponent{np^2U}\big)U^3}&\leqslant&
\Norm{\exponent{np^2U}\int_0^1(\Gamma_2U^2\Theta/
0.7)\exponent{\tau(\Gamma_2U^2\Theta/ 0.7)}\dd\tau
U^3}\nonumber\\
&\leqslant& C
np^3\norm{\exponent{0.5np^2U}U^5}\leqslant\frac{C}{p^2n\sqrt{n}}.
\label{s1}
\end{eqnarray}
Applying (\ref{NBVA}), (\ref{s1}) and Lemmas \ref{NBsmoothing} and
\ref{smoothing}, we obtain
 \begin{eqnarray*}
\lefteqn{\bigg\vert\norm{F-\NB(r,\qubar)}-\frac{\tilde
C_{TV}p}{\sqrt{n}}\bigg\vert\leqslant
\norm{F-\NB(r,\qubar)(\dirac+M_3)}+\bigg\vert\norm{\NB(r,\qubar)M_3}-\frac{\tilde
C_{TV}p}{\sqrt{n}}\bigg\vert}\hskip 1cm\\
&\leqslant&\frac{Cp}{n}+\norm{\NB(r,\qubar)(M_3+np^4U^3/3)}+\Ab{\frac{np^4}{3}\norm{\NB(r,\qubar)U^3}-\frac{\tilde
C_{TV}p}{\sqrt{n}}}\\
&\leqslant&
\frac{Cp^2}{\sqrt{n}}+\frac{np^4}{3}\norm{(\NB(r,\qubar)-\exponent{np^2U})U^3}+
\bigg\vert \frac{np^4}{3}\norm{\exponent{np^2U}U^3}-\frac{\tilde C_{TV}p}{\sqrt{n}}\bigg\vert\\
&\leqslant& \frac{Cp^2}{\sqrt{n}}+\frac{np^4}{3}
\bigg\vert\norm{\exponent{np^2U}U^3}-\frac{3\tilde
C_{TV}}{(np^2)^{3/2}}\bigg\vert\leqslant\frac{Cp^2}{\sqrt{n}}+\frac{C}{n}.
\end{eqnarray*}
 \qed


\textbf{Proof of Theorem \ref{k1k2th}.}  The direct consequence of
conditions  $(n-m+1)a(p)\geqslant 1$ and $ma(p)\leqslant 0.01$  are the following estimates
\[ (n-m+1)\geqslant 100m,\quad \tilde
N=\frac{(n-m+1)}{2m-1-m(m-1)/(n-m+1)}\geqslant\frac{100m}{2m}=50.\]
We have
\begin{eqnarray}
\pbar&=&\frac{(n-m+1)a(p)}{\tilde N}+\frac{(n-m+1)a(p)}{\tilde
N}\bigg(\frac{\tilde N}{N}-1\bigg)= \frac{(n-m+1)a(p)}{\tilde
N}\bigg(1+\frac{\epsilon}{\tilde
N-\epsilon}\bigg)\nonumber\\
&=&a(p)\bigg(2m-1-\frac{m(m-1)}{n-m+1}\bigg)\bigg(1+\frac{\epsilon}{\tilde
N-\epsilon}\bigg)\leqslant
a(p)\bigg(2m+\frac{m}{100}\bigg)\bigg(1+\frac{1}{49}\bigg)\nonumber\\
&\leqslant&2.05 a(p)m \leqslant 0.03.\label{barp1}
\end{eqnarray}

The sum $\tilde N$ has $n-m+1$ summands. After grouping, we get
$K$ $1$-dependent random variables containing $m$ initial summands
each, and (possibly) one additional variable, equal to the sum of
$\delta m$ initial summands. Here
\begin{equation}
K=\bigg\lfloor \frac{n-m+1}{m}\bigg\rfloor,\quad
\frac{n-m+1}{m}=K+\delta,\quad 0\leqslant\delta<1. \label{Kdelta}
\end{equation}
The analysis of the structure of new variables $X_j$ shows that,
for $j=1,\dots, K$
\begin{equation*}
X_j= \begin{cases} 1, & \mbox{with probability } ma(p), \\
0, & \mbox{with probability } 1-ma(p), \end{cases}\quad
X_{K+1}= \begin{cases} 1, & \mbox{with probability } \delta ma(p), \\
0, & \mbox{with probability } 1-\delta ma(p). \end{cases}
\end{equation*}
 Consequently,
 $\nu_2(j)=\nu_3(j)=\nu_4(j)=\wE_2^{+}(X_1,X_2)=\wE_2^{+}(X_1,X_2,X_3)=\wE_3(X_1,X_2)=0$.
 For calculation of $\Expect X_1X_2$, note that there are the
 following non-zero product events: a) the first summand of $X_1$ equals
 1 and any of the summands of $X_2$ equals 1 ($m$ variants); b) the second summand of $X_1$ equals
 1 and any of the summands of $X_2$, beginning from the second one,  equals 1 ($m-1$
 variant) and etc. Each event has the probability of occurrence $a^2(p)$.
 Therefore,
 \[\Expect
 X_1X_2=a^2(p)(m+(m-1)+(m-2)+\dots+1)=\frac{a(p)^2m(m+1)}{2}.\]
Similarly arguing we obtain the following relations for
$j=1,\dots,K$, ($j=2,\dots,K$ and $j=3,\dots,K$ if more variables
are involved) and $X_{K+1}$ (if $\delta>0$):
\begin{eqnarray}
\Expect X_j&=&ma(p),\quad\Expect
X_{j-1}X_j=\frac{m(m+1)a^2(p)}{2},\quad\wE(X_{j-1},X_j)=-\frac{m(m-1)a^2(p)}{2},\nonumber\\
\Expect X_{j-2}X_{j-1}X_j&=&\frac{m(m+1)(m+2)a^3(p)}{6},\quad
\wE(X_{j-2},X_{j-1},X_j)=\frac{a^3(p)m(m-1)(m-2)}{6},
\nonumber\\
\Expect X_{K+1}&=&\delta ma(p),\quad
 \Expect X_{K}X_{K+1}=\frac{\delta m(\delta m+1)a^2(p)}{2},\nonumber\\
\wE(X_{K},X_{K+1})&=&\frac{a^2(p)\delta m(\delta m+1-2m)}{2},\quad
\Expect X_{K-1}X_{K}X_{K+1}=\frac{\delta m(\delta m+1)(\delta m+2)a^3(p)}{6}\nonumber\\
\wE X_{K-1}X_{K}X_{K+1}&=&\frac{a^3(p)\delta
m(9m^2-9m+2)}{6}.\label{mix}
\end{eqnarray}
It is obvious, that $\Gamma_1=(n-m+1)a(p)$. Taking into account
(\ref{Kdelta}) and (\ref{mix}) we can calculate $\Gamma_2$:
\begin{eqnarray}
\Gamma_2&=&-\frac{1}{2}[Km^2a^2(p)+\delta^2m^2a^2(p)]-\frac{(K-1)m(m-1)a^2(p)}{2}+
\frac{\delta m a^2(p)(\delta m+1-2m)}{2}\nonumber\\
&=&-\frac{a^2(p)m}{2}[2m(K+\delta)-(K+\delta)-(m-1)]\nonumber\\
&=&-\frac{a^2(p)}{2}[(n-m+1)(2m-1)-m(m-1)]. \label{Gama2}
\end{eqnarray}
Similarly,
\begin{equation*}
\Gamma_3=\frac{a^3(p)}{6}[(n-m+1)(3m-1)(3m-2)-4m(2m-1)(m-1)].
\label{Gama3}
\end{equation*}
Making use of all the  formulas given above and noting that
$m\geqslant 2$, we to get the estimate
\begin{eqnarray*}
R_1&\leqslant& K(ma(p))^3+(\delta
ma(p))^3+3ma(p)[(K-2)m(m+1)a^2(p)/2+\delta
m(\delta m+1)a^2(p)/2]\nonumber\\
&&+C(K+\delta)m^3a^3(p)\leqslant Cm^3a^3(p)(K+\delta)\leqslant
C(n-m+1)m^2a^3(p). \label{R1a}\end{eqnarray*} Similarly,
\begin{equation*}
\wE^{+}(X_1,X_2,X_3,X_4)\leqslant C m^4a^4(p),\quad R_2\leqslant
C(n-m+1)m^3a^4(p). \label{R2a}
\end{equation*}
Using (\ref{Gama2}) and (\ref{barp1}), we get
\begin{eqnarray*}
\frac{N\pbar^3}{3}&=&\frac{\Gamma_1\pbar^2}{3}=
\frac{\Gamma_1}{3}\frac{4\Gamma_2^2}{\Gamma_1^2}\bigg(1+
\frac{\epsilon}{\tilde
N-\epsilon}\bigg)^2=\frac{4\Gamma_2^2}{3\Gamma_1}+\frac{4\Gamma_2^2}{3\Gamma_1}\frac{\epsilon}{\tilde
N-\epsilon}\bigg(2+\frac{\epsilon}{\tilde
N-\epsilon}\bigg)\nonumber\\
&=&  \frac{a^3(p)}{3}(n-m+1)(2m-1)^2 +\theta C m^3a^3(p).
\end{eqnarray*}
 Similarly,
\[\Gamma_3=\frac{a^3(p)}{6}(n-m+1)(3m-1)(3m-2)+\theta Cm^3a^3(p).\]
Therefore,
\begin{equation*}
\Gamma_3-\frac{N\pbar^3}{3}=A+C\theta m^3a^3(p). \label{skirt}
\end{equation*}
By Lemma \ref{a10}
\begin{equation}
m^3a^3(p)\norm{U^3\BI(N,\pbar)}\leqslant
C\frac{m^3a^3(p)}{(n-m+1)a(p)\sqrt{(n-m+1)a(p)}}\leqslant
C\frac{m^3a^2(p)}{n-m+1}.\label{liek}
\end{equation}
\noindent Next, we check  the conditions in (\ref{3ab}). Indeed,
we already noted that $\nu_2(j)=0$. Now
\begin{eqnarray*}
(K-1)\ab{\wE(X_1,X_2)}+\ab{\wE(X_{K-1},X_K)}&\leqslant&
\frac{Km(m-1)a^2(p)}{2}+\frac{\delta m 2ma^2}{2}
\leqslant(K+\delta)2m^2a^2(p)\\
&\leqslant&\frac{2ma^2}{n-m+1}=2ma\Gamma_1\leqslant 0.02\Gamma_1.
\end{eqnarray*}
It remains to apply Theorem \ref{Bintheorem} and (\ref{liek}).  \qed

\noindent \textbf{Proof of Theorem \ref{BINSHARP}.} We have
\begin{eqnarray*}
\lefteqn{\Ab{ \norm{\Ha-\BI(N,\bar{p})}- \tilde
C_{TV}\frac{a^{3/2}(p)m(m-1)}{2\sqrt{n-m+1}}} \leqslant
\norm{\Ha-\BI(N,\bar{p})\big(\dirac+ AU^3\big)}}\hskip
4cm\\
&&+\Norm{\BI(N,\pbar)U^3\Big(A-\frac{a^3(p)}{6}(n-m+1)m(m-1)\Big)}\\&&+
\frac{a^3(p)}{6}(n-m+1)m(m-1)\bigg\vert\norm{\BI(N,\pbar)U^3}-\frac{3\tilde
C_{TV}}{(N\pbar(1-\pbar))^{3/2}}\bigg\vert\\
&&+\frac{a^3(p)}{6}(n-m+1)m(m-1)\frac{3\tilde
C_{TV}}{(N\pbar)^{3/2}}\bigg\vert\frac{1}{(1-\pbar)^{3/2}}-1\bigg\vert.
\end{eqnarray*}
We easily check that
\[
\frac{1}{(1-\pbar)^{3/2}}-1=\frac{1-(1-\pbar)^3}{(1-\pbar)^{3/2}[1+(1-\pbar)^{3/2}]}=
\frac{\pbar[1+(1-\pbar)+(1-\pbar)^2]}{(1-\pbar)^{3/2}[1+(1-\pbar)^{3/2}]}=a(p)C(m)\theta.
\]
All that now remains is  to apply (\ref{BIVAk1k2}) and use Lemmas
\ref{a10} and \ref{sharpC}.  \qed

\vskip 0.5cm

\textbf{Acknowledgement.} We are grateful to the referees for useful remarks, which helped to improve the paper.


\end{document}